\documentclass[preprint,12pt]{elsarticle}
\usepackage{amsfonts,amsmath,amssymb}
\usepackage{tikz}

\newtheorem{thm}{Theorem}[section]
\newtheorem{lem}[thm]{Lemma}
\newtheorem{co}[thm]{Conjecture}
\newtheorem{prb}[thm]{Problem}
\newdefinition{df}{Definition}[section]
\newdefinition{rem}{Remark}[section]
\newdefinition{ex}{Example}[section]
\newproof{pf}{Proof}
\newproof{pot}{Proof of Theorem}
\numberwithin{equation}{section}

\journal{
}

\begin{document}

\begin{frontmatter}

\title{
Criteria of irreducibility of the Koopman representations for the group ${\rm GL}_0(2\infty,{\mathbb R})$
}
\author{A.V.~Kosyak
}
\ead{kosyak02@gmail.com}
\address{Max-Planck-Institut f\"ur Mathematik, Vivatsgasse 7, D-53111 Bonn, Germany}
\address{Institute of Mathematics, Ukrainian National Academy of Sciences,
3 Tereshchenkivs'ka Str., Kyiv, 01601, Ukraine}

\begin{abstract}
Our aim is to find the irreducibility criteria for the Koopman representation, when the
group acts on some space with a measure (Conjecture~\ref{co.G-Ism-art}).
Some general necessary conditions of the irreducibility of this representation are established.
In the particular case of the group ${\rm
GL}_0(2\infty,{\mathbb R})$ $= \varinjlim_{n}{\rm
GL}(2n-1,{\mathbb R})$, the inductive limit of the general linear
groups we prove that these conditions are also the necessary ones.
The corresponding measure is infinite tensor products of
one-dimensional arbitrary Gaussian non-centered measures. The
corresponding $G$-space $X_m$ is a subspace of the space ${\rm
Mat}(2\infty,{\mathbb R})$ of infinite in both directions real
matrices. In fact, $X_m$  is a collection of $m$ infinite in both
directions  rows.
This result was announced in \cite{Kos04}. We give the proof only for $m\leq 2$. The general
case will be studied later.
\end{abstract}


\end{frontmatter}
\newpage
\tableofcontents
\section{ Introduction }


%
\subsection{Description of the dual for locally compact groups}
The main problem in the representation theory for
a locally compact group  $G$ is to find  the set of {\it all unitary
irreducible representations} of $G$ up to unitary equivalence and to decompose
reducible representations into a direct sum or direct integral of irreducible.  This set is called the {\it unitary dual} of $G$ and is denoted by $\hat
G$. For many locally compact groups this problem has been solved,
but for some particular cases it remains open, for example, for the group ${\rm
SO}(p,q)$. To find the dual for locally compact groups
$G$, one can use {\it regular, quasiregular or induced
representations}. In the case of locally compact groups all these
constructions are based on the existence of the invariant {\it Haar
measure} on the initial group $G$ or some $G$-quasi-invariant
measure on the corresponding homogeneous space $H\setminus G$,
where $H$ is a closed subgroup of $G$ or on  some general $G$-space $X$.

\subsection{Regular, quasiregular and induced representations for infinite-dimensional groups}
It is well known that there is no
general method to describe $\hat G$ for
infinite-dimensional groups $G$.
Our aim is to start the development of the harmonic analysis on
infinite-dimensional groups.

In the previous articles we have  generalized the notions of the regular, quasiregular and induced
representations for infinite-dimensional groups by constructing
$G$-quasi-invariant measures on  suitable completions of the
corresponding objects (groups, homogeneous spaces and $G$-spaces).
In addition, we study the irreducibility of the constructed
representations in the framework of the Ismagilov conjecture (see
\ref{co.Ism-art}).

In this article we consider the case when the infinite-dimensi\-onal group $G$, the {\it inductive limit of the general linear groups}, acts on the space of $m$ infinite rows equipped with the Gaussian measure.
We establish  the {\it criteria of irreducibility}
of constructed representations (see Theorem~\ref{5.t.irr}) in terms of the corresponding measure
 and express some general conjectures dealing with the irreducibility.
These conjectures are natural generalization of the {\it Ismagilov conjecture} (see Conjecture~\ref{co.G-Ism-art}).

Recall some previous constructions. {\it Regular representations} for
infinite-dimensio\-nal groups were defined and studied in \cite{
Kos90,Kos92,Kos94}.
Due to the result of A.Weil \cite{Weil53}, there is no invariant measure on  non
locally compact groups. Therefore, to construct an analogue of a regular
representation of an infinite-dimensional group $G$ we can, for
example,  construct a $G$-quasi-invariant measure on a suitable
completion $\tilde G$ of the  initial group $G$.
The regular representation of an
infinite-dimensional group can be irreducible, which never happens
for a locally compact group, except for the trivial one!

To define
a {\it quasiregular representation} we should construct a
$G$-quasi-invariant measure on a suitable completion $\tilde
H\setminus \tilde G$ of the homogeneous space $H\setminus G$
\cite{Kos02.2,Kos02.3,Kos03}.

To construct the induced
representation for infinite-dimensional groups we need to extend
by continuity the representation of the subgroup $H$ to the
corresponding completion $\tilde H$. The general construction of the {\it induced
representations} and the beginning of the
{\it orbit methods}  for infinite-dimensional group of upper triangular matrices were
done in  \cite{KosJFA14}.

To construct the regular representation for an infinite-dimensi\-onal
group $G$, first we should find some larger topological group $\widetilde{G}$
and a measure $\mu$ on $\widetilde{G}$ such that $G$ is a dense
subgroup in $\widetilde{G},$ and $\mu^{R_t}\sim\mu$ for all $t\in
G,$ (or $\mu^{L_t}\sim\mu$ for all $t\in G$), here $\sim$ means {\it equivalence}. The right and left
representations $T^{R,\mu},T^{L,\mu}:G\rightarrow U(L^2(\tilde
G,\mu))$ are naturally defined in the Hilbert space $L^2(\tilde
G,\mu)$ by the following formulas:
$$
(T^{R,\mu}_tf)(x)=(d\mu(xt)/d\mu(x))^{1/2}f(xt),
$$
$$
(T^{L,\mu}_sf)(x)=(d\mu(s^{-1}x)/d\mu(x))^{1/2}f(s^{-1}x).
$$
The right regular representation of infinite-dimensional
groups can be irreducible if no left actions are {\it admissible} for
the measure $\mu$, i.e., when $\mu^{L_t}\perp \mu$ for all $t\in
G\backslash{\{e\}}$. In this case a von Neumann algebra
${\mathfrak A}^{T^{L,\mu}}$ generated by the left regular
representation $T^{L,\mu}$  is trivial.  More precisely:
\begin{co}
[Ismagilov, 1985] \label{co.Ism-art} The right regular
representation \vskip -0.2cm
$$T^{R,\mu}:G\rightarrow U(L^2(\tilde G,\mu))$$
is irreducible if and only if \par 1) $\mu^{L_t}\perp
\mu\,\,\forall t\in G\backslash{\{e\}},\,\,$ (where $\perp$ stands
for singular),
\par
2) the measure $\mu$ is $G$-ergodic.
\end{co}
This conjecture was verified for a lot of particular cases. In the general case, it is an open problem.
In the case of a finite field ${\mathbb F}_p$ we need some additional conditions for the irreducibility  \cite{Kos10}.

\subsection{Koopman representation}

Let $\alpha:G\rightarrow {\rm Aut}(X)$ be a measurable action of a
group $G$  on a measurable space $(X,\mu)$ with
$G$-quasi-invariant measure $\mu$, i.e, $\mu^{\alpha_t}\sim\mu$
for all $t\in G$. With these date one can associate the
representation $ \pi^{\alpha ,\mu,X}:G\rightarrow U(L^2(X,d\mu)),
$ by the following formula:
\begin{equation}
\label{Rep(G,X)-pi}
(\pi^{\alpha,\mu,X}_tf)(x)=(d\mu(\alpha_{t^{-1}}(x))/d\mu(x))^{1/2}f(\alpha_{t^{-1}}(x)),\quad
f\in L^2(X,\mu).
\end{equation}
In the case of an invariant measure this representation
 called {\it Koopman's representation}, see  \cite{Koo31}.
\index{represenation!Koopman's}
\index{subgroup!centarlizer}
%
%
We would like to solve the following problems:
\begin{prb}
\label{p.Koop1-?}
Find criteria of irreducibility of the representation $\pi^{\alpha ,\mu,X}$ defined by (\ref{Rep(G,X)-pi}).
\end{prb}
\begin{prb}
\label{p.Koop2-?}
Find the description of the commutant of the von Neumann algebra generated by representation $\pi^{\alpha ,\mu,X}$ when
representation is reducible.
\end{prb}

To study properties of the Koopman representation, in particular, the
irreducibility, {\it we need some conjectures to describe the
commutant of the von  Neumann} algebras generated by this
representation.
The Schur--Weyl duality and the Dixmier commutation theorem
below give us a very good hint for such a conjecture, see
Conjecture~\ref{co.Com} in a general context.

\subsection{Schur--Weyl duality}
\index{duality!Schur--Weyl}  {\it Schur--Weyl duality} \cite{Schur01,Schur27,Weyl39} is a typical
situation in representation theory involving two kinds of symmetry
that determine each other.

From \cite{Yuan12}: ``If $V$ is a finite-dimensional complex vector space, then the symmetric group $S_n$
naturally acts on the tensor power $V^{\otimes n}$ by permuting the factors. This action of $S_n$ commutes
with the action of ${\rm GL}(V)$, so all permutations $\sigma : V^{\otimes n} \to V^{\otimes n}$ are
morphisms of ${\rm GL}(V)$-representations. This defines a morphism $\mathbb{C}[S_n] \to
\text{End}_{\text{GL}(V)}(V^{\otimes n})$, and a natural question to ask is whether this map is surjective.

Part of Schur--Weyl duality asserts that the answer is yes. The {\it double commutant theorem}
plays an important role in the proof and also highlights an important corollary, namely that
$V^{\otimes n}$ admits a canonical decomposition
$$
V^{\otimes n} = \bigoplus_{\lambda} V_{\lambda} \otimes S_{\lambda}
$$
where $\lambda$ runs over partitions,
$V_{\lambda}$ are some irreducible representations of ${\rm GL}(V)$,
and $S_{\lambda}$
are the {\it Specht modules}, which describe
all irreducible representations of $S_n$. This gives a fundamental
relationship between the representation theories of the general linear and
symmetric groups; in particular, the assignment $V \mapsto V_{\lambda}$ can be upgraded to a functor
called a {\it Schur functor}, generalizing the construction of the exterior and symmetric products.''
\index{module!Specht}
\index{functor!Schur}
\index{theorem!double commutant}

Let ${\rm dim}V=m$ then ${\rm GL}(V)={\rm GL}(m,\mathbb C)$.
The abstract form of the Schur--Weyl duality asserts that two algebras of operators on the tensor space generated by the
actions of ${\rm GL}(m,\mathbb C)$ and $S_n$ are the full mutual centralizers in the algebra of the endomorphisms
${\displaystyle \mathrm {End} _{\mathbb {C} }(\mathbb {C}^m\otimes \mathbb {C}^m\otimes \cdots \otimes \mathbb {C}^m)}$.

Denote by $\alpha$ and $\beta$ the corresponding actions of $S_n$
and ${\rm GL}(m,\mathbb C)$ in the group of all automorphisms ${\rm
Aut}({\mathbb C}^m\otimes  {\mathbb C}^m\otimes \cdots
\otimes{\mathbb C}^m)$:
$$
\alpha:S_n\to {\rm Aut}(X),\quad \beta:{\rm GL}(m,\mathbb C)\to
{\rm Aut}(X).
$$
Let $M'$ be the commutant of the subset $M$ in the von Neumann
algebra $B(H)$ of all bounded operators in a Hilbert space $H$:
\begin{equation}
\label{M^'}
M'=\{B\in B(H)\mid[B,a]=0\,\, \forall a\in M\}\,\,\text{where}\,\,\,[B,a]=Ba-aB.
\end{equation}
Set $M_1=(\alpha(S_n))''$ and $M_2=(\beta({\rm GL}(m,{\mathbb
C})))''$ then the Schur--Weyl duality states that
$M_1'=M_2\quad \text{hence,}\quad M_2'=M_1.$

In \cite{TsiVer14} the authors extend the classical Schur--Weyl duality between representations of the groups  $\text{SL}(m,\mathbb C)$ and  $S_n$  to the case of  $\text{SL}(m,\mathbb C)$  and the infinite symmetric group  $S_\infty$.
In \cite{PenSty11} the authors extend Weyl results to the classical infinite-dimensional locally finite algebras
${\mathfrak gl}_\infty,\,\,{\mathfrak sl}_\infty,\,\,{\mathfrak sp}_\infty,\,\,{\mathfrak so}_\infty$.

\subsection{The Dixmier commutation theorem, locally compact groups}
Let $G$ be a locally compact group and let $h$ be
the right invariant Haar measure on $G$, i.e., $h^{R_t}=h$ for all $t\in G$. Consider the left $L$ and the right $R$ action of the group $G$ on itself:
$$
R_t(x)=xt^{-1},\,\,L_s(x)=sx,\,\,x,t,s\in G.
$$
The right and the left regular representations of the group $G$
are defined in the Hilbert space $L^2(G,h)$ by
$$
(\rho_tf)(x)=f(xt),\quad
(\lambda_sf)(x)=\big(dh(s^{-1}x)/dh(x)\big)^{-1/2}f(s^{-1}x),\,\,f\in
L^2(G,h),
$$
where $dh(s^{-1}x)/dh(x)$ is the Radon-Nikodim derivative.
\begin{thm} [Dixmier's commutation theorem \cite{Dix69C}]
\label{t.Dix[,]}
The commutant of the von-Neumann algebra generated by the right
regular representation is generated by the left regular
representation. More precisely, let $\rho,\lambda :G \rightarrow
U(L^2(G,h))$ be the right and the left regular representations of
the group $G$, and let  ${\mathfrak A}^\rho=(\rho_t\mid t\in G)''$ and
${\mathfrak A}^\lambda=(\lambda_s\mid s\in G)''$ be the
corresponding von Neumann algebras. Then
\begin{equation}
\label{Dix:A'=B}
({\mathfrak A}^\rho)'={\mathfrak A}^\lambda\quad  \text{and}\quad
 ({\mathfrak A}^\lambda)'={\mathfrak A}^\rho.
\end{equation}
\end{thm}

\subsection{$G$-action and  irreducibility of the Koopman representation}
%
In both  examples we have two commuting actions of the group $G_1$ and $G_2$ on the same space $X$.
Let $Z_{G}(H)$ be a {\it centralizer} of the subgroup $H$ in the group $G$:
$$
Z_{G}(H)=\{g\in G\mid \{g,a\}=e\,\,\forall a\in H\},
$$
where $\{g,a\}=gag^{-1}a^{-1}$.
In the first example, we have two commuting actions $\alpha$ and $\beta$ of the groups
$G_1=S_n$ and $G_2={\rm GL}(n,\mathbb C)$ on the space $X$ such that
$Z_{{\rm Aut}(X)}(\alpha(G_1))\supseteq\beta(G_2)$. In the second example, we have two commuting actions
$R$ and $L$ of the same group $G$ in the space $X=G$. In this case we have $\{R(G),L(G)\}=e$ or
$Z_{{\rm Aut}(G)}(R(G))\supseteq L(G)$. In the general case, if we have only one group $G$ acting via $\alpha$ on the space
$X$, the second group should be the {\it centralizer} of the group $\alpha(G)$ in the group ${\rm Aut}(X)$, i.e.,
it is natural to consider $G_2=Z_{{\rm Aut}(X)}(\alpha(G))$.

Come back to the Koopman representation (\ref{Rep(G,X)-pi}).
Consider the centralizer $Z_{{\rm Aut}(X)}(\alpha(G))$   of the
subgroup $\alpha(G)=\{\alpha_t\mid t\in G\}$ in the group ${\rm
Aut}(X)$ and its subgroup $G_2$ defined as follows:
$$
G_2:=Z^\mu_{{\rm Aut}(X)}(\alpha(G)):=\big\{g\in Z_{{\rm
Aut}(X)}(\alpha(G))\mid \mu^g\sim\mu\big\}.
$$
Define  the representation $T$ of the group $G_2$ as follows:
\begin{equation}
\label{Rep(T)}
(T_gf)(x)=(d\mu(gx)/d\mu(x))^{1/2}f(gx).
\end{equation}
Consider two von Neumann algebras
$$
{\mathfrak A}^{\pi}(G) =(\pi_t\mid t\in G)'',\quad {\mathfrak
A}^{T}(G_2)=(T_g\mid g\in G_2)''.
$$
The conditions 1) and 2) below are necessary conditions of the irreducibility of the representation $\pi^{\alpha ,\mu,X}$.
It would be interesting to know when they are sufficient, i.e., when the following conjecture is true
\begin{co} [Kosyak, \cite{Kos94,Kos02.3}]
\label{co.G-Ism-art}
The representation
$$\pi^{\alpha ,\mu,X}:G\rightarrow U(L^2(X,\mu))
$$
is irreducible if and only if
\par 1) $\mu^g\perp \mu\,\,\forall
g\in Z_{{\rm Aut}(X)}(\alpha(G))\backslash{\{e\}},\,\,$
\par 2) the measure $\mu$ is $G$-ergodic.
\end{co}
Recall that a measure $\mu$ is $G$-ergodic if
$f(\alpha_t(x))\!=\!f(x)\,\,\mu$ a.e. for all $t\in G$ implies
$f(x)\!=\!const\,\,\mu$ a.e.(almost everywhere) for all functions $f\!\in\! L^1(X,\mu)$. 
\begin{co}
\label{co.Com} The commutant of the von Neumann algebra generated
by representation $\pi$ (\ref{Rep(G,X)-pi}) of the group $G$ coincides  with the von Neumann algebra generated by the representation $T$ (\ref{Rep(T)}) of the subgroup $G_2$ in the centralizer
$Z_{{\rm Aut}(X)}(\alpha(G))$:
$$
({\mathfrak A}^{\pi}(G))' ={\mathfrak A}^{T}(G_2).
$$
\end{co}
For a lot of particular cases Conjecture \ref{co.Com} holds,
but in general it fails. Below we give several example
for which Conjecture~\ref{co.Com} fails.

\subsection{Counterexample to  Conjecture~\ref{co.Com}}
\label{ch753}
\subsubsection{Case $X=S_{n-1}\setminus S_n$}
\begin{ex}
\label{ex.S(3)/S(2)}
Consider the group $S_n$ acting on the homogeneous space $X=S_{n-1}\setminus S_n$.
For corresponding right quasiregular representation of $S_n$  in $L^2(X)$ Conjecture~\ref{co.Com} fails.
\end{ex}
\begin{pf}
To simplify details set $n=3$. For general $n$ the proof is the same.
Let $\sigma_1,\sigma_2$ be two generators of the group $S_3$:
\begin{equation}
\label{S_3}S_3=\Big(\sigma_1,\sigma_2\mid \sigma_1^2=e,\sigma_2^2=e,\,
\sigma_1\sigma_2\sigma_1=\sigma_2\sigma_1\sigma_2\Big).
\end{equation}
Let the group $S_2$ is generated by $\sigma_1$, then the space $X$ consists of three classes $x_0=\{e, \sigma_1\},\,\,x_1=\{\sigma_2,\sigma_1\sigma_2\},\,\,
x_3=\{\sigma_2\sigma_1,\sigma_1\sigma_2\sigma_1\}$. The right action of $S_3$ on the space $X$ is as follows:
\begin{eqnarray*}
x_0\sigma_1=x_0,\quad x_1\sigma_1=x_2,\quad  x_2\sigma_1=x_1,\\
x_0\sigma_2=x_1,\quad x_1\sigma_2=x_0,\quad  x_2\sigma_2=x_2.
\end{eqnarray*}
Therefore, in $L^2(X)$ the corresponding representations for $T_{\sigma_1}$ and $T_{\sigma_2}$ are as follows:
$$
T_{\sigma_1}=
\left(\begin{smallmatrix}
 1&0&0\\
 0&0&1\\
 0&1&0
 \end{smallmatrix}\right),
\quad
T_{\sigma_2}=
\left(\begin{smallmatrix}
 0&1&0\\
 1&0&0\\
 0&0&1
  \end{smallmatrix}\right).
$$
The representation $T$ is reducible, since the vector $e_0+e_1+e_2$ is invariant. It splits into one-dimensional and two- dimensional irreducible representations.
 But the group $S_3$ acts on $X$ by permutations so, its centralizer is trivial.
\qed\end{pf}
\subsubsection{Case $X=O(3)\backslash O(3)$}
\begin{ex}
\label{ex.SO/SO}
Consider the group $O(3)$ acting on the homegeneous space $O(3)\backslash O(3)$ $\simeq S^2$. The centralizer of $O(3)$ in the group of all automorphisms ${\rm Aut}(S^2)$ consists of two elements $I$ and $-I$ by Lemma~\ref{l.Z(SO(3))/)} but the representation of $O(3)$ in $L^2(X)$ is an infinite direct sum of irreducible representations
generated by eigenvectors
of the Laplace operator on $S^2$, see \cite[Chapter I,\S 3]{Hel84}. Therefore, Conjecture~\ref{co.Com} fails.
\end{ex}
\subsubsection{Centralizer of $SO(2k+1)$}
Let $n\geq0$, and $SO(n)$ be the group of all real orthogonal
$n\times n$-matrices with determinant $1$. This group effectively and
transitively acts $n-1$-dimensional sphere $S^{n-1}$, and so it
can be regarded as a subgroup of the group $\mathcal{H}(S^{n-1})$
of all homeomorphisms of $S^{n-1}$.

Let $I$ be the unit matrix.
Then $-I$ is an ``antipodal'' map, that is $-I(p) = -p$ for all $p\in S^{n-1}$.
Evidently, $I$ and $-I$ commute with all elements from $SO(n)$, and so $\{\pm I\}$ belongs to the centralizer of $SO(n)$ in $O(n)$.


\begin{lem}
\label{l.Z(SO(3))/)}
Suppose $n=2k+1$ is odd.
Then the group $\{\pm I\}$ is the centralizer of $SO(2k+1)$ in all the group $\mathcal{H}(S^{2k})$.
\end{lem}
\begin{pf} (given by S.~Maximenko.)
Suppose $h\in\mathcal{H}(S^{n-1})$ commutes with all matrices $A \in SO(n)$, that is $h\circ A(x) = A \circ h(x)$ for all $x\in S^{n-1}$.
We should prove that then $h = \pm I$.

First we claim that $h(x) \in \{\pm x\}$ for each $x\in S^{n-1}$.
Indeed, since $n$ is odd, for each $x\in S^{n-1}$ there exists $A\in SO(n)$ such that $\{\pm x\}$ is the set of all fixed points for $A$.
Hence
\[h(x) = h \circ A(x) = A \circ h(x),\]
that is $h(x)$ is a fixed point for $A$, and so $h(x) = \pm x$.

Now, suppose $h(x) = \varepsilon x$ for some $\varepsilon = \pm 1$.
We claim that then $h = \varepsilon I$.
Let $F = \{ x \in S^{n-1} \mid h(x) = \varepsilon x\}$ be the set of points where $h$ coincides with $\varepsilon I$.
We will show that $F$ is a non-empty open-closed subset of $S^{n-1}$, which will imply that $F$ coincides with all of $S^{n-1}$.

As shown above $x\in F$, so $F\not=\varnothing$.
Moreover, as $h$ and $-I$ are continuous, $F$ is closed.
It remains to show that $F$ is open.
Let $U$ be a small neighbourhood of $x$ such that $U \cap -U = \varnothing$, that is $U$ does not contain antipodal pairs.
Since $h$ is continuous and $h(x) = \varepsilon x \in \varepsilon U$, there exists a neighbourhood $V$ of $x$ such that $h(V)=\varepsilon U$.
Then for each $y\in V$ we have that $h(y) \in \{\pm y\} \cap \varepsilon U = \varepsilon y$.
In other words, $h=\varepsilon I$ on $V$, and so $V\subset F$.
This proves that $F=S^{n-1}$.
\qed\end{pf}

\section{Representations of the inductive limit of the general linear groups
${\rm GL}_0(2\infty,{\mathbb R})$} \label{sec5.1}
\subsection{Finite-dimensional case}
Consider the space
$
X_{m,n}=\Big\{x=\sum_{1\leq k\leq m}\sum_{-n\leq r\leq n} x_{kr}E_{kr},\,\,x_{kr}\in
{\mathbb R}\Big\},
$
with the measure   (see (\ref{5.mu^m}))
$
\mu_{(b,a)}^{m,n}(x)=\otimes_{k=1}^m\otimes_{-n\leq r\leq n}\mu_{(b_{kr},a_{kr})}(x_{kr}).$
On the space $X_{m,n}$ acts two groups $\text{GL}(m,\mathbb R)$ from the left and $\text{GL}(2n+1,\mathbb R)$ from the right and these actions commute. Therefore, two von Neumann algebras ${\mathfrak A}_1$  and ${\mathfrak A}_{2,n} $ in the Hilbert space $L^2(X_{m,n},\mu_{(b,a)}^{m,n})$  generated respectively by the left and the right actions of the corresponding groups have the property that ${\mathfrak A}_1'\subseteq {\mathfrak A}_{2,n}$. We study what happens when $n\to\infty$. As the limit we obtain some unitary representation of the group  ${\rm GL}_0(2\infty,{\mathbb R})=\varinjlim_{n,i^s}{\rm GL}(2n-1,{\mathbb R})$ (see below). In generic case, this representation is reducible, namely, if there exists a non trivial element $s\in \text{GL}(m,\mathbb R)$ such the the left action is admissible for the measure $\mu_{(b,a)}^m$, i.e., $(\mu_{(b,a)}^m)^{L_s}\sim \mu_{(b,a)}^m$. But when no non-trivial left actions are admissible, i.e.,
when $(\mu_{(b,a)}^m)^{L_s}\perp \mu_{(b,a)}^m$ for all $s\in \text{GL}(m,\mathbb R)\backslash \{e\}$  we prove that this representation is irreducible Theorem~\ref{5.t.irr}. Here, as in the case of the regular \cite{Kos90,Kos92} and   quasiregular \cite{Kos02.2,Kos02.3} representations of the group $B^{\mathbb N}_0$ we obtain the remarkable fact that  the irreducible representations can be obtained as the inductive limit of reducible representations!
\subsection{Infinite-dimensional case}
Let us denote by ${\rm Mat}(2\infty,{\mathbb R})$ the space of all
real  matrices infinite in both directions:
\begin{equation}
\label{5.Mat(2inf,R)} {\rm Mat}(2\infty,{\mathbb R})=
\Big\{x=\sum_{ k,n\in{\mathbb Z}}x_{kn}E_{kn},\,\,x_{kn}\in
{\mathbb R}\Big\},
\end{equation}
where $E_{kn},\,\,k,n\in {\mathbb Z}$ are infinite matrix unities.

\index{${\rm GL}_0(2\infty,{\mathbb R})=\varinjlim_{n,i^s}{\rm
GL}(2n+1,{\mathbb R})$}
The group ${\rm GL}_0(2\infty,{\mathbb R})=\varinjlim_{n,i^s}{\rm
GL}(2n-1,{\mathbb R})$ is defined as the inductive limit of the
general linear groups $G_n={\rm GL}(2n-1,{\mathbb R})$ with
respect to the symmetric embedding $i^s$ (\ref{N.i^s}):
\begin{equation}
\label{N.i^s} {\rm GL}(2n-1,{\mathbb R})\ni x\mapsto
i^s_{n+1}(x)=x+E_{-(n+1),-(n+1)}+E_{n+1,n+1}\in {\rm GL}(2n+1,{\mathbb R}).
\end{equation}
 We consider a $G$-space $X_m,\,\,m\in {\mathbb N}$  as the
following subspace of the space ${\rm Mat}(2\infty,{\mathbb R})$:
\begin{equation}
\label{5.X_m} X_m=\Big\{x\in {\rm Mat}(2\infty,{\mathbb
R})\,\vert\, x=\sum_{k=1}^m\sum_{n\in{\mathbb
Z}}x_{kn}E_{kn}\Big\}.
\end{equation}
The  group ${\rm GL}_0(2\infty,{\mathbb R})$ acts from the right
on the space $X_m.$ Namely, the right action of the group $G={\rm
GL}_0(2\infty,{\mathbb R})$ is correctly defined on the space
$X_m$ by the formula $R_t(x)=xt^{-1},\,\,t\in G,\,\,x\in X_m$. We
define a Gaussian noncentered product measure
$\mu^m=\mu_{(b,a)}^m$ on the space $X_m:$
\begin{equation}
\label{5.mu^m}
\mu_{(b,a)}^m(x)=\otimes_{k=1}^m\otimes_{n\in{\mathbb
Z}}\mu_{(b_{kn},a_{kn})}(x_{kn}),
\end{equation}
 where
$$
d\mu_{(b_{kn},a_{kn})}(x_{kn})=(b_{kn}/\pi)^{1/2}
\exp(-b_{kn}(x_{kn}-a_{kn})^2)dx_{kn}
$$
and $b=(b_{kn})_{k,n},\,\,b_{kn}>0,\,a=(a_{kn})_{k,n},\,a_{kn}\in
{\mathbb R},\,1\leq k\leq m,\,n\in {\mathbb Z}.$ Define the
representation $T^{R,\mu,m}$ of the group ${\rm
GL}_0(2\infty,{\mathbb R})$ in the space $L^2(X_m,\mu^m_{(b,a)})$
by the formula:
$$
(T^{R,\mu,m}_tf)(x)=\big(d\mu_{(b,a)}^m(xt)/d\mu_{(b,a)}^m(x)\big)^{1/2}f(xt),\,\,f\in L^2(X_m,\mu^m_{(b,a)}).
$$
Obviously, the centralizer $Z_{{\rm Aut}(X_m)}(\alpha(G))\subset
{\rm Aut}(X_m)$ contains  the group $L({\rm GL}(m,{\mathbb R}))$,
i.e., the image of the group ${\rm GL}(m,{\mathbb R})$ with
respect to the left action $L:{\rm GL}(m,{\mathbb R})\rightarrow
{\rm Aut}(X_m),\,L_s(x)\!=\!sx,\,s\in {\rm GL}(m,{\mathbb
R}),\,x\in X_m.$ We prove the following theorem  for $m\leq 2$.
\begin{thm}
\label{5.t.irr} The representation $T^{R,\mu,m}\!:\!{\rm
GL}_0(2\infty,{\mathbb R})\!\rightarrow\! U(L^2(X_m,\mu^m_{(b,a)}))$ is
irreducible if and only if $(\mu^m_{(b,a)})^{L_{s}}\perp
\mu^m_{(b,a)}\,\, \forall s\in {\rm GL}(m,{\mathbb
R})\backslash\{e\}.$
\end{thm}
\begin{rem}
Any Gaussian product-measure $\mu_{(b,a)}^m$ on $X_m$ is ${\rm
GL}_0(2\infty,{\mathbb R})$-right-ergodic \cite[\S 3, Corollary 1]{ShFDT67}. For non-product-measures
this is not true in general.
\end{rem}
To study the condition $(\mu^m_{(b,a)})^{L_{t}}\perp
\mu^m_{(b,a)}$ for $t\in {\rm GL}(m,{\mathbb R})$ set
\begin{equation}
\label{5.X_n(t)}
t=(t_{rs})_{r,s=1}^m\in{\rm GL}(m,{\mathbb R}), \,\,\,B_n={\rm
diag}(b_{1n},b_{2n},...,b_{mn}),\,\,\,X_n(t)=B_n^{1/2}tB_n^{-1/2}.
\end{equation}
Let  $M^{i_1i_2...i_r}_{j_1j_2...j_r}(t)$ be the minors of
the matrix $t$ with $i_1,i_2,...,i_r$ rows and $j_1,j_2,...,j_r$
columns, $1\leq r\leq m.$ Let $\delta_{rs}$ be the Kronecker
symbols.
\begin{lem}
\label{perp1} For the measures $\mu_{(b,a)}^m,\,m\in{\mathbb N}$
the relation
 $(\mu_{(b,a)}^m)^{L_t}\perp\mu_{(b,a)}^m\,$ $ \forall t\in
 {\rm GL}(m,{\mathbb R})\backslash\{e\}$ holds if and only if
$$
\prod_{n\in {\mathbb Z}}\frac{1}{2^m\vert{\rm det}\,\,t\vert} {\rm
det}\left(I+X_n^*(t)X_n(t)\right)+ \sum_{n\in {\mathbb
Z}}\sum_{r=1}^m
b_{rn}\left(\sum_{s=1}^m(t_{rs}-\delta_{rs})a_{sn}\right)^2
=\infty,
$$
$$
\quad{\rm where}\quad {\rm det}\left(I+X_n^*(t)X_n(t)\right)=
$$
$$
 1+\sum_{r=1}^{m} \sum_{1\leq i_1<i_2<...<i_r\leq m;1\leq
j_1<j_2<...<j_r\leq m}
\left(M^{i_1i_2...i_r}_{j_1j_2...j_r}(X_n(t))\right)^2.
$$
\end{lem}
This lemma will be proved in Section~\ref{ch5.6}.
\begin{rem}
\label{r.irr-idea-G} {\rm The idea of the proof of the
irreducibility.} Let us denote by ${\mathfrak A}^m$ the  von
Neumann algebra generated by the representation $T^{R,\mu,m}:$ $
{\mathfrak A^m}=(T^{R,\mu,m}_t\mid t\in G)''. $ For
$\alpha\!=\!(\alpha_k)\!\in\!\{0,1\}^m$ define the von Neumann algebra
$L^\infty_\alpha(X_m,\mu^m)$ as follows:
$$
L^\infty_\alpha(X_m,\mu^m)\!=\!\big(\exp(itB^\alpha_{kn})\mid
1\leq k\leq m,\,\,t\in {\mathbb R},\,\,n\in {\mathbb
Z}\big)'',\,\,
$$
where $\text{ }\,\,
B^\alpha_{kn}\!=\!\left\{\begin{array}{ccc}
 x_{kn},&\text{if}&\alpha_k=0\\
 D_{kn},&\text{if}&\alpha_k=1
\end{array}\right..$

The proof of the {\rm irreducibility is based on  three
facts}:\par 1) using the orthogonality condition
$(\mu^m)^{L_{t}}\perp \mu^m$ for all $t\in {\rm GL}(m,{\mathbb
R})\backslash{\{e\}}$ we can approximate by generators
$A_{kn}=A_{kn}^{R,m}=\frac{d}{dt}T^{R,\mu,m}_{I+tE_{kn}}\vert_{t=0}$
the set of operators $(B^\alpha_{kn})_{k=1}^m,\,n\!\in\!{\mathbb
Z}$ for some $\alpha\!\in\!\{0,1\}^m$ depending on the measure
$\mu^m$,\par 2) the subalgebra $L^\infty_\alpha(X_m,\mu^m)\subset
{\mathfrak A}^m$ is a maximal abelian subalgebra in ${\mathfrak
A}^m$,\par
3)  the measure $\mu^m$ is $G$-ergodic. \\
Here the generators $A_{kn}^{R,m}$ are given by the formulas:
$$
A_{kn}^{R,m}\!=\!\sum_{r=1}^{m}x_{rk}D_{rn},\quad k,n \in
{\mathbb Z},\quad\text{where}\quad D_{kn}=\partial/\partial x_{kn}-b_{kn}(x_{kn}-a_{kn}).
$$
\end{rem}
\begin{rem}
\label{r.key-lemma2}
The fact that conditions $(\mu^m)^{L_{t}}\!\perp\! \mu^m$ for all $t\in {\rm GL}(m,{\mathbb
R})\backslash{\{e\}}$ implies the possibility of the approximation of $x_{kn}$ and $D_{kn}$ is based on some
completely independent statement about {\rm the properties of projections of two infinite vectors} $f=(f_k)_{k\in{\mathbb N}}$ and
$g=(g_k)_{k\in{\mathbb N}}$ such that $f,\,g,\,f+sg\not\in l_2$ for all $s\in {\mathbb R}$ (Lemma~\ref{l.min=proj}). This lemma is a key part of the proof of the irreducibility of the representation.
\end{rem}
\begin{rem}
\label{r.key-lemma3}
Similarly, for the ``nilpotent group'' $B_0^{\mathbb N}$ and the infinite  product of arbitrary  Gaussian measures on ${\mathbb R}^m$ (see \cite{KosAlb06J}) the proof of the irreducibility is based on another {\rm completely independent statement} namely, {\rm Hadamard -- Fischer's inequality}, see Lemma~\ref{l.det(det)1}.

\begin{lem}[Hadamard -- Fischer's inequality
\cite{HornJon89}, \cite{HornJon91}
]
\label{l.det(det)1}  For any
positive definite matrix $C\in{\rm Mat}(m,{\mathbb
R}),\,\,m\in{\mathbb N}$ and  any two subsets $\alpha$ and $\beta$
with $\emptyset\subseteq\alpha,\,\,\beta\subseteq\{1,...,m\}$ the
following inequality holds:
\begin{equation}
\label{det(det)1}
 \left|\begin{array}{cc}
M(\alpha)&M(\alpha\bigcap\beta)\\
M(\alpha\bigcup\beta)&M(\beta)
\end{array}\right|=
 \left|\begin{array}{cc}
A(\hat{\alpha})&A(\hat{\alpha}\bigcup\hat{\beta})\\
A(\hat{\alpha}\bigcap\hat{\beta})&A(\hat{\beta})
\end{array}\right|
\geq 0
\end{equation}
where $M(\alpha)=M^\alpha_\alpha(C),\,\,
A(\alpha)=A^\alpha_\alpha(C)$ and
$\hat{\alpha}=\{1,...,m\}\setminus\alpha.$
\end{lem}
For details see \cite[p.573]{HornJon89}, \cite[Chapter
2.5, problem 36]{HornJon91}.

The conditions of orthogonality $\mu^{L_{t}}\perp \mu$ with respect to the left action of the group $B(m,{\mathbb R})$ on $X^m$ were expressed as the divergence of some series, $S^L_{kn}(\mu)=\infty,\,\,1\leq k<n\leq m$. Conditions of the approximation of the variables $x_{kn}$ by combinations of generators $A_{pq}$ were expressed in terms of the divergence of another series $\Sigma_{kn}$. The proof of the fact that conditions $S^L_{kn}(\mu)=\infty,\,\,1\leq k<n\leq m$ imply conditions
$\Sigma_{kn}=\infty,\,\,1\leq k<n\leq m$ is based on the
Hadamard -- Fischer's inequality.
\end{rem}
\index{inequality!Hadamard -- Ficher's}
\section{The  proof of the irreducibility}
\subsection{The cases $m=1$}
\label{sec5.2}As before, let us denote by $\langle f_n\mid n\in {\mathbb
N}\rangle$ the closure of the linear space generated by the set of
vectors $(f_n)_{n\in{\mathbb N}}$ in a Hilbert space $H.$ We shall write $\mu_{(b,a)}=\mu_{(b,a)}^1$.

{\bf In the case $m=1$} the generators $A_{kn}^{R,1}$ have the
form
$$
A_{kn}^{R,1} =x_{1k}D_{1n},\qquad k,n \in {\mathbb Z}.
$$
The following lemmas are proved in \cite{KosAlb05J}
\begin{lem}
\label{5.1} The following  three conditions are equivalent:\par
(i)\,$\,\,(\mu_{(b,a)})^{L_t}\perp\mu_{(b,a)}$ for all $t\in
GL(1,{\mathbb R})\setminus\{e\}$, \par
(ii)\,\,\,$(\mu_{(b,a)})^{L_{-E_{11}}}\perp\mu_{(b,a)}$,\par
(iii)\,\,\,\,$S^L_{11}(\mu)=4\sum_{n\in{\mathbb Z}}{b_{1n}}a_{1n}^2=\infty. $
\end{lem}
\begin{lem}
\label{5.2}
 For $k,m\in {\mathbb Z}$ we have
$$
x_{1k}x_{1m}{\bf 1}\in\langle A_{kn}^{R,1}A_{mn}^{R,1}{\bf
1}=x_{1k}x_{1m}D_{1n}^2{\bf 1} \mid n\in {\mathbb Z}\rangle.
$$
\end{lem}
\begin{lem}
\label{5.3}
 For any $k\in {\mathbb Z}$ we have
$$
x_{1k}{\bf 1}\in\langle x_{1k}x_{1n}{\bf 1}\mid n\in {\mathbb Z}\rangle
\Leftrightarrow S^L_{11}(\mu)=\infty.
$$
\end{lem}
So, operators $x_{1k},\,\,k\in {\mathbb Z}$  are affiliated (see
\cite{Dix69W}) with the von Neumann algebra ${\mathfrak A}^1$
(notation $x_{1k}\,\,\eta\,\, {\mathfrak A}^1$) which completes
the proof of the  irreducibility for $m=1.$
\section{The  proof of the irreducibility in the cases $m=2$ }
\label{sec.5.3}
{\bf In the case $m=2$} the generators
$A_{kn}:=A_{kn}^{R,2}:=\frac{d}{dt}T^{R,\mu,2}_{I+tE_{kn}}\mid_{t=0}$
have the form:
$$
A_{kn}=x_{1k}D_{1n}+x_{2k}D_{2n},\qquad k,n \in {\mathbb Z}.
$$
\begin{lem}
\label{perp2}
 Three following conditions (i)--(iii) are equivalent for
the measure $\mu=\mu_{(b,a)}^2${\rm :}
\par (i) $\mu^{L_t}\perp\mu$ for all
$t\in GL(2,{\mathbb R})\setminus\{e\},$ where $L_t(x)=tx,\,\,x\in
X_2;$
$$
(ii) \left\{\begin{array}{cll}
 (a)&\mu^{L_{\exp(tE_{12})}}\perp\mu,    &\forall t\in {\mathbb R}\backslash\{0\},\\
 (b)&\mu^{L_{\exp(tE_{21})}}\perp\mu,    &\forall t \in {\mathbb R}\backslash\{0\},\\
 (c)&\mu^{L_{\exp(tE_{12})P_1}}\perp\mu, &\forall t\in {\mathbb R},               \\
 (d)&\mu^{L_{\exp(tE_{21})P_2}}\perp\mu, &\forall t\in {\mathbb R},               \\
 (e)&\mu^{L_{\tau_{-}(\phi,s)}}\perp\mu, & \forall \tau_{-}(\phi,s) \in GL(2,{\mathbb R})\setminus\{e\},
\end{array}\right.
$$
$$
(iii) \left\{\begin{array}{cll}
(a)& S^L_{12}(\mu)=\infty,&\\
(b)&S^L_{21}(\mu)=\infty,&\\
(c)&S^{L,-}_{12}(\mu,t)=\infty,&\forall t\in {\mathbb R},\\
(d)&S^{L,-}_{21}(\mu,t)=\infty,&\forall t\in {\mathbb R}, \\
(e)&
\Sigma^-_{12}(\tau_{-}(\phi,s))=\infty,&
\forall s>0,\phi\in[0,2\pi),
\end{array}\right.
$$
 where
\begin{equation}
\label{S^L_(kn)}
S^L_{kn}(\mu)=\sum_{m\in {\mathbb Z}}
\frac{b_{km}}{2}\left(\frac{1}{2b_{nm}}+a_{nm}^2\right),\quad k\not=n,
\end{equation}
\begin{equation}
\label{S^(L,-)_(kn)}
S^{L,-}_{kn}(\mu,t)=\frac{t^2}{4}\sum_{m\in{\mathbb
Z}}\frac{b_{km}}{b_{nm}}+ \sum_{m\in{\mathbb
Z}}\frac{b_{km}}{2}(-2a_{km}+ta_{nm})^2,
\end{equation}
$$
\Sigma^-_{12}(\tau_{-}(\phi,s))=\sin^2\phi\Sigma_1(s)+\Sigma^-_2(\tau_{-}(\phi,s)),
$$
\begin{equation}
\label{sigma1(s)}
\Sigma_1(s):=\!\sum_{n\in {\mathbb Z}}
\Big(s^2\sqrt{\frac{b_{1n}}{b_{2n}}}\!-\!s^{-2}\sqrt{\frac{b_{2n}}{b_{1n}}}\Big)^2,
\end{equation}
\begin{equation}
\label{sigma2(s)} \Sigma^-_2(\tau_{-}(\phi,s))\!:=\!\sum_{n\in
{\mathbb
Z}}\big(4\sin^2\frac{\phi}{2}b_{1n}\!+\!4\cos^2\frac{\phi}{2}s^{-4}b_{2n}\big)
\big(\sin\frac{\phi}{2}a_{1n}\!-\!s^2\cos\frac{\phi}{2}a_{2n}\big)^2,
\end{equation}
$$
\exp(tE_{12})=I+tE_{12}= \left(\begin{array}{cc}
 1&t\\
 0&1
\end{array}\right),\quad
\exp(tE_{21})=I+tE_{21}=\left(\begin{array}{cc}
 1&0\\
 t&1
\end{array}\right),
$$
$$
\exp(tE_{12})P_1= \left(\begin{array}{cc}
 -1&t\\
  0&1
\end{array}\right),\quad
\exp(tE_{21})P_2= \left(\begin{array}{cc}
 1&0\\
 t&-1
\end{array}\right),
$$
$$
\tau_{-}(\phi,s)= \left(
\begin{array}{cc}
\cos\phi&s^2\sin\phi\\
s^{-2}\sin\phi&-\cos\phi
\end{array}
\right)\,\,
\text{and}\,\,
P_1= \left(\begin{array}{cc}
 -1&0\\
  0&1
\end{array}\right),\,\,
P_2=
 \left(\begin{array}{cc}
 1& 0\\
 0&-1
\end{array}\right).
$$
Moreover, $(ii)(\sharp)\Leftrightarrow(iii)(\sharp)$ for
$\sharp=a,b,c,d,e.$
\end{lem}
\begin{rem}
\label{5.5}
 We observe that
$$
\tau_{-}(\phi,s)\!=\!\left(\!
\begin{array}{cc}
\cos\phi&s^2\sin\phi\\
s^{-2}\sin\phi&-\cos\phi
\end{array}
\!\right)\!=\!\left(\!\begin{array}{cc}
s&0\\
0&s^{-1}
\end{array}\!\right)\!
\left(\!
\begin{array}{cc}
\cos\phi&-\sin\phi\\
\sin\phi&\cos\phi
\end{array}
\!\right)\! \left(\!
\begin{array}{cc}
s^{-1}&0\\
0     &s
\end{array}
\!\right)\!P_2.
$$
\end{rem}
\begin{rem}
\label{r.5.5}
We note \cite[Chapter V ,\S 8 Problems, 2, p. 147]{Knap86} that every element of   ${\rm SL}(2,{\mathbb R})$ is
conjugate to at least one matrix of the form
$$
\left(\begin{array}{cc}
 a&0\\
  0&a^{-1}
\end{array}\right),\,\,a\not=0,\,\,\left(\begin{array}{cc}
 1&t\\
  0&1
\end{array}\right),\,\,\left(\begin{array}{cc}
 -1&t\\
  0&-1
\end{array}\right),\,\,
\left(\begin{array}{cc}
\cos\phi&\sin\phi\\
-\sin\phi&\cos\phi
\end{array}
\right).
$$
\end{rem}
%
\begin{rem}
\label{perp2-1} The three following conditions are equivalent:
\begin{eqnarray*}
(i)&\quad\mu^{L_{\tau_-(\phi,s)}}\perp\mu,&\phi\in[0,2\pi),\,\,s>0,\\
(ii)&\Sigma_{12}^-(\tau_{-}(\phi,\!s))\!=\!\sin^2\phi\Sigma_1(s)\!+\!\Sigma^-_2(\tau_{-}(\phi,s))\!=\!
\infty,&\phi\in[0,2\pi),\,\,s>0,\\
(iii)&\quad\Sigma_1(s)+\Sigma_2(C_1,C_2)=\infty,&s\!>\!0,\,(C_1,C_2)\!\in\!{\mathbb R}^2\!\setminus\!\{0\},
\end{eqnarray*}
where $\Sigma_1(s)$ is defined by (\ref{sigma1(s)}) and
$$
\Sigma_2(C_1,C_2):=
\sum_{n\in {\mathbb Z}}(C_1^2b_{1n}+C_2^2b_{2n})(C_1a_{1n}+C_2a_{2n})^2.
$$
\end{rem}
\begin{pf}
In Section~\ref{ch5.6} we shall show that $(i)\Leftrightarrow (ii)$ (see
(\ref{tau(phi,s)-})), i.e., that
$$
\mu^{L_{\tau_-(\phi,s)}}\perp\mu \Leftrightarrow
\Sigma^-_{12}(\tau_{-}(\phi,s))= \sin^2\phi\Sigma_1(s)+\Sigma^-_2(\tau_{-}(\phi,s))= \infty.
$$
To prove $(ii)\Leftrightarrow (iii)$ set
$$
\sin\frac{\psi}{2}\!=\!\sin\frac{\phi}{2}(\sin^2\frac{\phi}{2}\!+\!s^4\cos^2\frac{\phi}{2})^{-1/2},\,\,\,
 \cos\frac{\psi}{2}\!=\!s^2\cos\frac{\phi}{2}(\sin^2\frac{\phi}{2}\!+\!s^4\cos^2\frac{\phi}{2})^{-1/2}
$$
then using (\ref{sigma2(s)}) we get
$$
\Sigma^-_2(\tau_{-}(\phi,s))\!:=
(\sin^2\frac{\psi}{2}+s^4\cos^2\frac{\psi}{2})^24
\!\sum_{n\in
{\mathbb
Z}}\big(\sin^2\frac{\psi}{2}b_{1n}\!+\!s^{-8}\cos^2\frac{\psi}{2}b_{2n}\big)\times
$$
$$
\big(\sin\frac{\psi}{2}a_{1n}\!-\!\cos\frac{\psi}{2}a_{2n}\big)^2\sim\Sigma_2(\psi):=\sum_{n\in
{\mathbb
Z}}\big(\sin^2\frac{\psi}{2}b_{1n}+\cos^2\frac{\psi}{2}b_{2n}\big)\times
$$
$$
\big(\sin\frac{\psi}{2}a_{1n}-\cos\frac{\psi}{2}a_{2n}\big)^2\!=\!\Sigma(C_1,C_2)\!=\!\sum_{n\in {\mathbb Z}}(C_1^2b_{1n}+C_2^2b_{2n})(C_1a_{1n}+C_2a_{2n})^2.
$$
\qed\end{pf}
%
\subsection{Some orthogonality problem in  measure theory}
Our aim now is to find the minimal set of conditions of the orthogonality
$\mu^{L_t}\perp\mu$ for all
$t\in GL(2,{\mathbb R})\setminus\{e\}$. To be more precise, consider more general situation.

Let $\alpha:G\rightarrow {\rm
Aut}(X)$  be a measurable action of a group $G$  on a
measurable space $(X,{\mathfrak B},\mu)$ with the following property:
 $\mu^{\alpha_t}\perp\mu$ for all $t\in G\setminus\{e\}$. Consider a subset $G^\perp(\mu)$ in the group $G$
having the following property:
\begin{equation}
\label{min-perp}
\text{if}\quad \mu^{\alpha_t}\perp\mu\,\,\forall t\in G^\perp(\mu)\quad\text{then}\quad \mu^{\alpha_t}\perp\mu\,\,\forall t\in G\setminus\{e\}.
\end{equation}
{\bf Problem.} Find a minimal subset $G^\perp_0(\mu)$ having the
property (\ref{min-perp}).
\begin{ex}
\label{nil(G)} Consider the nilpotent group $B(m,{\mathbb R})$ of
upper triangular real $m\times m$ matrices with units on the
diagonal acting on the space $X_m$ with the Gaussian product
measure $\mu=\mu^m_{(b,a)}$, where $X_m$ and $\mu$ are defined as
follows (see details in  \cite{Kos02.2,Kos02.3}):
$$
X_m=\{I+\sum_{1\leq k\leq m}\sum_{k<n}x_{kn}E_{kn}\},\quad \mu^m=\otimes_{1\leq k\leq m}\otimes_{k<n}\mu_{(b_{kn},a_{kn})}.
$$
\end{ex}
Using results form \cite{Kos02.2}  and \cite{Kos02.3}
we conclude that the three following conditions are equivalent:
\begin{eqnarray*}
(i)&\mu^{L_t}\perp\mu&\forall t\in B(m,{\mathbb R}) \backslash\{e\},\\
(ii)&\mu^{L_{\exp(tE_{kn})}}\perp\mu       &\forall t\in {\mathbb R}\backslash\{0\},\quad       1\leq k<n\leq m,\\
(iii)& S_{kn}^L(\mu)=\infty                 &1\leq k<n\leq m,
\end{eqnarray*}
where $S_{kn}^L(\mu)$ is defined by (\ref{S^L_(kn)})
 $$ S_{kn}^L(\mu)
=\sum_{r=n+1}^\infty\frac{b_{kr}}{2}\Big(\frac{1}{2b_{nr}}+a_{nr}^2\Big).
$$
In fact, it is sufficient to fix a nontrivial point $t_{kn}\not=0$
on  any subgroup $\exp(tE_{kn})=I+tE_{kn},\,\,t\in {\mathbb R},\,\,1\leq k<n\leq m$. In this case
the subset $G_0^\perp(\mu)$ is discrete and consists of
$m(m-1)/2$ points:
$$
G_0^\perp(\mu, t)\!=\!\Big(I+t_{kn}E_{kn}\mid t_{kn}\in {\mathbb R}\backslash\{0\}, 1\leq k<n\leq m\,\,\Big),
$$
where $t\!=\!(t_{kn})_{kn}\in ({\mathbb R}\backslash\{0\})^{m(m-1)/2}$.
For $ t_1\not = t_2\in ({\mathbb R}\backslash\{0\})^{m(m-1)/2}$ we get two distinct minimal subsets
$G_0^\perp(\mu^m,t_1)$ and $G_0^\perp(\mu^m,t_2)$.
\begin{ex}
\label{solvable(G)}
 Consider the solvable group $Bor(m,{\mathbb R})$ of
upper triangular real $m\times m$ matrices with nonzero elements
on the diagonal acting on the space $X_m$ with the Gaussian
product measure $\mu=\mu^m_{(b,a)}$, where $X_m$ and $\mu$ are
defined as follows (see details in \cite{KosAlb05J})
$$
X_m=\{x=\sum_{1\leq k\leq m}\sum_{k\leq n}x_{kn}E_{kn}\},\quad
\mu^m_{(b,a)}=\otimes_{1\leq k\leq m}\otimes_{k\leq n
}\mu_{(b_{kn},a_{kn})}.
$$
Using \cite[Theorem 5]{KosAlb05J}
we conclude that the following three conditions are equivalent:
\begin{eqnarray*}
(i)&\mu^{L_t}\perp\mu, &\forall t\in Bor(m,{\mathbb R}) \backslash\{e\},\\
(ii)&\mu^{L_{\exp(tE_{kn})}}\perp\mu \,\,      \forall t\in {\mathbb
R}\backslash\{0\},&1\leq k<n\leq m,\\
   &\mu^{L_{\exp(tE_{kn})P_k}}\perp\mu \,\,\forall t\in {\mathbb
R},&1\leq k<n\leq m,\\
(iii)& S_{kn}^L(\mu)=\infty ,\quad S^{L,-}_{kn}(\mu,t)=\infty, &1\leq k<n\leq m,
\end{eqnarray*}
 where $S^{L,-}_{kn}(\mu,t)$ is defined by (\ref{S^(L,-)_(kn)}).  As before, it is sufficient to fix a nontrivial point $t_{kn}\not=0$
on  any subgroup $\exp(tE_{kn})=I+tE_{kn},\,\,t\in {\mathbb R}$. But on the curves $\exp(tE_{kn})P_k$ we can not omit  any point
$t\in {\mathbb R}$. Finally, a minimal subset depending on the choice of $t=(t_{kn})_{kn}\in ({\mathbb R}\backslash\{0\})^{m(m-1)/2}$
can be chosen as follows:
$$
G_0^\perp(\mu, t)\!=\!\Big(\exp(t_{kn}E_{kn})=I+t_{kn}E_{kn}\mid t_{kn}\in {\mathbb R}\backslash\{0\}, 1\leq k<n\leq m\,\,\Big)\bigcup
$$
$$
\Big(\exp(tE_{kn})P_k\mid \,\,\forall t\in {\mathbb R}, 1\leq k<n\leq m\,\,\Big)
$$
where $P_k=I-2E_{kk}$. For example, for $m=2$ we get $P_1={\rm diag}(-1,1)$ and $P_2={\rm diag}(1,-1)$.
\end{ex}
\begin{ex}
\label{GL} In the case of the group $G={\rm GL}(2,{\mathbb R})$
acting on the space  $X_2$ defined by (\ref{5.X_m}) with the
measure $\mu_{(b,a)}^2$ defined by (\ref{5.mu^m}) using
Lemmas~\ref{perp2}, we conclude that the
description of the set $G_0^\perp(\mu_{(b,a)}^2)$ is as follows:
\begin{equation*}
G_0^\perp(\mu_{(b,a)}^2,t_{12},t_{21})=\Big(\exp(t_{12}E_{12}),\,\exp(t_{21}E_{21})\mid t_{12},t_{21}\in {\mathbb R}\backslash\{0\}\Big)
\bigcup
\end{equation*}
\begin{equation*}
                 \Big(\exp(tE_{12})P_1,\,\exp(tE_{21})P_2\mid \forall t\in {\mathbb R}\Big)\bigcup
                 \Big(\tau_-(\phi,s)\mid \forall s>0,\,\,\phi\in[0,2\pi)\Big).
\end{equation*}
\end{ex}
\begin{rem}
\label{r.g^2=e}
We note that except the one-parameter groups $E_{kn}(t)=I+tE_{kn},\,t\in {\mathbb R}$ all other element from the set
$G_0^\perp(\mu)$ for $G={\rm GL}(2,{\mathbb R})$ are of order 2, i.e., if $g\in \{\exp(tE_{kn})P_k,\,\,\tau_-(\phi,s)\}$ then $g^2=e$.
\end{rem}
\subsection{Approximation of $x_{kn}$ and $D_{kn}$}
We will formulate several  lemmas, which will be  useful for
approximation of the independent variables $x_{kn}$ and operators
$D_{kn}$ by combinations of the generators $A_{kn}$. For short, we
shall write $A_{kn}$ instead of $A_{kn}^{R,2}.$

In what follows we use the following notation
for $f,g\in{\mathbb R}^m$
\begin{equation}
\label{D(f,g)}
\Delta(f,g)=\frac{\Gamma(f)+\Gamma(f,g)}{\Gamma(g)+1}.
\end{equation}
\begin{lem}
\label{x1x1}
  For any  $k,t\in {\mathbb Z}$ one has
$$
x_{1n}x_{1t}\in\langle A_{nk}A_{tk}{\bf 1}\mid k\in{\mathbb
Z}\rangle \,\,\Leftrightarrow\,\,
\lim_m\Sigma_{1,m}(x,x)=\infty,
$$
where $\Sigma_{1,m}(x,x)=\Delta(f_m^1,g_m^1)$ and
\begin{equation}
\label{f^1,g^1+}
f_m^1=\Big(\frac{b_{1k}}{\sqrt{b_{1k}^2+2b_{1k}b_{2k}}}\Big)_{k=-m}^m,\quad g_m^1=\Big(\frac{b_{2k}}
{\sqrt{b_{1k}^2+2b_{1k}b_{2k}}}\Big)_{k=-m}^m.
\end{equation}
\end{lem}
\begin{lem}
\label{x2x2} For any  $k,t\in {\mathbb Z}$ we have
$$
x_{2k}x_{2t}\in\langle A_{kn}A_{tn}{\bf 1}\mid n\in {\mathbb
Z}\rangle\,\,\Leftrightarrow\,\,\lim_m\Sigma_{2,m}(x,x)=\infty,
$$
where $\Sigma_{2,m}(x,x)=\Delta(f_m^2,g_m^2)$ and
\begin{equation}
\label{f^2,g^2+}
f_m^2=\Big(\frac{b_{2k}}{\sqrt{b_{2k}^2+2b_{1k}b_{2k}}}\Big)_{k=-m}^m,\quad g_m^2=\Big(\frac{b_{1k}}{\sqrt{b_{2k}^2+2b_{1k}b_{2k}}}\Big)_
{k=-m}^m.
\end{equation}
\end{lem}
\begin{rem}
\label{x1x1>}
We say that two series $\sum_na_n$ and $\sum_nb_n$ with positive $a_n,\,b_n$ are {\it equivalent} if they are simultaneously convergent or divergent. In this case we shall use the notaions $\sum_na_n\sim \sum_nb_n$.
Using the obvious equivalence of the following two series with positive $a_n$ and $b_n$
\begin{equation}
\label{sim1}
\sum_{n\in{\mathbb N}}\frac{a_n}{a_n+b_n}\sim \sum_{n\in{\mathbb N}}\frac{a_n}{b_n}
\end{equation}
we have the following estimation (we set $\Sigma^{12}=\sum_{k\in{\mathbb Z}}\frac{b_{1k}}{b_{2k}}$ and
$\Sigma^{21}=\sum_{k\in{\mathbb Z}}\frac{b_{2k}}{b_{1k}}$)
\begin{eqnarray*}
\Vert f^1\Vert^2\!=\!\sum_{k\in{\mathbb Z}}\frac{b_{1k}^2}{b_{1k}^2+2b_{1k}b_{2k}}\!\sim\!\sum_{k\in{\mathbb Z}}\frac{b_{1k}}{2b_{2k}}
\!=\!\frac{\Sigma^{12}}{2},\\
\Vert f^2\Vert^2\!=\!\sum_{k\in{\mathbb Z}}\frac{b_{2k}^2}{b_{2k}^2+2b_{1k}b_{2k}}\!\sim\!\sum_{k\in{\mathbb Z}}\frac{b_{2k}}{2b_{1k}}
\!=\!\frac{\Sigma^{21}}{2},\\
\Vert g^1\Vert^2\!=\!\sum_{k\in{\mathbb Z}}\frac{b_{2k}^2}{b_{1k}^2+2b_{1k}b_{2k}}\! <\! \sum_{k\in{\mathbb Z}}\frac{b_{2k}}{2b_{1k}}
\!=\!\frac{\Sigma^{21}}{2},\\
\Vert g^1\Vert^2\!=\!\sum_{k\in{\mathbb Z}}\frac{b_{1k}^2}{b_{2k}^2+2b_{1k}b_{2k}}\! <\!\sum_{k\in{\mathbb Z}}\frac{b_{1k}}{2b_{2k}}\!=\!
\frac{\Sigma^{12}}{2},
\end{eqnarray*}
we conclude that $\lim_{m}\Sigma_{1,m}(x,x)\!=\!\infty$ if
$$
\quad\lim_m\Sigma_{1,m}'(x,x):=\!\lim_m\Big(\sum_{k=-m}^m \frac{b_{1k}}{b_{2k}}\Big)
\Big(\sum_{k=-m}^m \frac{b_{2k}}{b_{1k}}\Big)^{-1}\!=\!\Sigma^{12}/\Sigma^{21}\!=\!\infty
$$
and $\lim_{m}\Sigma_{2,m}(x,x)\!=\!\infty$ if
$$
\lim_{m}\Sigma_{2,m}'(x,x):=\!\lim_{m}\Big(\sum_{k=-m}^m \frac{b_{2k}}{b_{1k}}\Big)
\Big(\sum_{k=-m}^m \frac{b_{1k}}{b_{2k}}\Big)^{-1}\!\!=\!\Sigma^{21}/\Sigma^{12}\!=\!\infty.
$$
\end{rem}
\begin{lem}
\label{d1}
 For any  $n\in {\mathbb Z}$ we have
$$
D_{1n}{\bf 1}\in\langle A_{kn}{\bf 1}\mid k\in {\mathbb Z}\rangle
\quad\Leftrightarrow\quad \lim_m\Sigma_{1,m}(D)=\infty,
$$
where $\Sigma_{1,m}(D)=\Delta(f_m,g_m)$ and
\begin{equation}
\label{f_m,g_m=+}
f_m\!=\!\Big(a_{1k}\Big(\frac{1}{2b_{1k}}\!+\!\frac{1}{2b_{2k}}\Big)^{-1/2}\Big)_{k=-m}^m,\,\,\,
g_m\!=\!\Big(a_{2k}\Big(\frac{1}{2b_{1k}}\!+\!\frac{1}{2b_{2k}}\Big)^{-1/2}\Big)_{k=-m}^m.
\end{equation}
\end{lem}

\begin{lem}
\label{d2} Set $\Sigma_{2,m}(D)=\Delta(g_m,f_m)$. For any  $n\in {\mathbb Z}$ we get
$$
D_{2n}{\bf 1}\in\langle A_{kn}{\bf 1}\mid k\in {\mathbb Z}\rangle
\quad\Leftrightarrow\quad \lim_m\Sigma_{2,m}(D)=\infty.
$$
\end{lem}

\begin{lem}
\label{x1,dA}
 For any  $k\in {\mathbb Z}$ we get
 $$
 x_{1k}{\bf 1}\in\langle D_{1n}A_{kn}{\bf 1}\mid n\in{\mathbb Z}\rangle
 \quad\Leftrightarrow\quad
\sum_{n\in{\mathbb Z}}\frac{b_{1n}}{b_{2n}}=\infty.
$$
\end{lem}
\begin{lem}
\label{x2,dA}
 For any  $k\in {\mathbb Z}$ we have
 $$
 x_{2k}{\bf 1}\in\langle D_{2n}A_{kn}{\bf 1}\mid n\in{\mathbb Z}\rangle
 \quad\Leftrightarrow\quad
\sum_{n\in{\mathbb Z}}\frac{b_{2n}}{b_{1n}}=\infty.
$$
\end{lem}

Our aim now is to show that some of the expressions
$\Sigma_{1,m}(x,x),\,\, \Sigma_{2,m}(x,x)$ and $\Sigma_{1,m}(D)$,
$\Sigma_{2,m}(D)$ tend to infinity if $\mu^{L_t}\perp\mu$ for all
$t\in GL(2,{\mathbb R})\setminus\{e\}$.

Let $\Gamma(f_1,f_2,\dots,f_n)$ be the Gramm determinant and
$\gamma(f_1,f_2,\dots,f_n)$ be the Gramm matrix of $n$ vectors
$f_1,f_2,\dots,f_n$ in a Hilbert space (see \cite{Gan58}).
The following lemma is trivial and well known but
we need  exact formulas.
\begin{lem}
\label{dist1} Let $f_1,f_2$ be two vectors in a Hilbert space. The distance
$\delta\langle f_2,f_1\rangle$ of the vector $f_2$ from the line
$\langle f_1\rangle$ generated by $f_1$ is given by the following
formula:
\begin{equation}
\label{dist-(f,g)}
\delta^2\langle f_2,f_1\rangle=\Vert f_2-\frac{(f_2,f_1)}{(f_1,f_1)}f_1\Vert^2=\frac{\Gamma(f_1,f_2)}{\Gamma(f_1)}.
\end{equation}
\end{lem}
\begin{pf} Obviously,
$\delta^2\langle f_2,f_1\rangle\!=\!\Vert f_2\!-\!f_0\Vert^2$ where $f_0\!=\!C_1f_1$ such that \\$(f_2\!-\!f_0,f_1)\!=\!0$. We have
$$
0=(f_2-f_0,f_1)=(f_2,f_1)-C_1(f_1,f_1)\quad\text{hence},\quad C_1=\frac{(f_2,f_1)}{(f_1,f_1)}.
$$
Finally,
\begin{eqnarray*}
\delta^2\langle f_2,f_1\rangle=
\Vert f_2-f_0\Vert^2=\Vert f_2-C_1f_1\Vert^2=
(f_2,f_2)-2C_1(f_2,f_1)+C_1^2(f_1,f_1)=\\
(f_2,f_2)
-\frac{2(f_2,f_1)(f_2,f_1)}{(f_1,f_1)}+\frac{(f_2,f_1)^2}{(f_1,f_1)^2}(f_1,f_1)=\\
\frac{(f_2,f_2)(f_1,f_1)-(f_2,f_1)(f_1,f_2)}{(f_1,f_1)}=
\frac{\Gamma(f_1,f_2)}{\Gamma(f_1)}.
\end{eqnarray*}
\qed\end{pf}
%
\begin{lem}
\label{l.min=proj0}
Let $f=(f_k)_{k\in{\mathbb N}}$ and $g=(g_k)_{k\in{\mathbb N}}$ be two real vectors  such that
$\Vert f\Vert^2=\infty$  where $\Vert f\Vert^2=\sum_kf_k^2$.  Denote by $f_{(n)}$,
$g_{(n)}\in {\mathbb R}^n$ their projections to the subspace
${\mathbb R}^n$, i.e., $f_{(n)}=(f_k)_{k=1}^n,\quad  g_{(n)}=(g_k)_{k=1}^n$ and set
\begin{equation}
\label{Delta-to-infty}
\Delta(f_{(n)},g_{(n)})=\frac{\Gamma(f_{(n)})+\Gamma(f_{(n)},g_{(n)})}{\Gamma(g_{(n)})+1}
\quad\text{then}\quad \lim_{n\to\infty}\Delta(f_{(n)},g_{(n)})=\infty
\end{equation}
in the following cases:
\begin{eqnarray*}
(a)&\Vert g\Vert^2<\infty,\\
(b)&\Vert g\Vert^2=\infty,\quad\text{and}\quad  \lim_{n\to\infty}\frac{\Vert f_{(n)}\Vert}{\Vert g_{(n)}\Vert}=\infty,\\
(c)&\Vert f\Vert^2=\Vert g\Vert^2=\Vert f+s g\Vert^2=\infty,\quad\text{for all}\quad
     s\in {\mathbb R}\setminus\{0\}.
\end{eqnarray*}
\end{lem}
\begin{pf} Obviously $\lim_{n\to\infty}\Delta(f_{(n)},g_{(n)})=\infty$ if conditions (a)  or (b) hold.
The implication $ (c)\Rightarrow (\ref{Delta-to-infty})$ is based on the following lemma.
\qed\end{pf}
\begin{lem}
\label{l.min=proj}
Let $f=(f_k)_{k\in{\mathbb N}}$ and $g=(g_k)_{k\in{\mathbb N}}$ be two real vectors  such that
\begin{equation}
\label{norm=infty}
\Vert f\Vert^2=\Vert g\Vert^2=\Vert C_1 f+C_2 g\Vert^2=\infty,\quad\text{for all}\quad
(C_1,C_2)\in {\mathbb R}^2\setminus\{0\},
\end{equation}
\begin{equation}
\label{final.}
\text{then}\quad
\lim_{n\to\infty}\frac{\Gamma(f_{(n)},g_{(n)})}{\Gamma(g_{(n)})}=\infty\quad\text{and}\quad
\lim_{n\to\infty}\frac{\Gamma(f_{(n)},g_{(n)})}{\Gamma(f_{(n)})}=\infty.
\end{equation}
\end{lem}
\begin{pf}
Assume that $\frac{\Vert f_{(n)}\Vert}{\Vert g_{(n)}\Vert}\leq C_1,\,\, \forall n\in{\mathbb N}$. The case $\frac{\Vert f_{(n)}\Vert}{\Vert g_{(n)}\Vert}\geq C_1$
is similar.  In this case
$
\frac{\Gamma(f_{(n)},g_{(n)})}{\Gamma(g_{(n)})}\leq C_1^2 \frac{\Gamma(f_{(n)},g_{(n)})}{\Gamma(f_{(n)})}
$
therefore, to prove (\ref{final.}) it is sufficient to prove that
$\lim_{n\to\infty}\frac{\Gamma(f_{(n)},g_{(n)})}{\Gamma(g_{(n)})}=\infty$.
Let us suppose the opposite, i.e., that  for all $n\in {\mathbb N}$ holds
\begin{equation}
\label{G(fg):G(g)<C.1}
\quad\quad\frac{\Gamma(f_{(n)},g_{(n)})}{\Gamma(g_{(n)})}\leq C.
\end{equation}
Set
$t_n=\frac{\Vert f_{(n)}\Vert}{\Vert g_{(n)}\Vert}$ then by the inequality $\frac{\Vert f_{(n)}\Vert}{\Vert g_{(n)}\Vert}\leq C_1$
there exists a subsequence $t_{n_k}$ such that the limit exists
$$
\lim_{k\to \infty}t_{n_k}=t_0\in [0,C_1].
$$
Let $\alpha_n$ be an angle between two vectors $f_{(n)},g_{(n)}\in{\mathbb R}^n$.
Since $\frac{\Gamma(f,g)}{\Gamma(g)}$ is the square of the distance of the vector $f$ from the line generated by $g$ by Lemma~\ref{dist1},
we have
\begin{equation}
\label{alpha(n)-to-0}
\frac{\Gamma(f_{(n)},g_{(n)})}
{\Gamma(g_{(n)})}=\Vert f_{(n)}\Vert^2\sin^2\alpha_n\leq C,\quad
\text{therefore}\quad \alpha_n\sim\Vert f_n\Vert^{-1}\to 0.
\end{equation}
 For $k,n\in{\mathbb N}$ set $M(k,n)=\left|\begin{array}{cc}
 f_k&f_n\\
  g_k&g_n
\end{array}\right|$, then by the  Lagrange identity for $f_{(m)}=(f_k)_{k=1}^m,\,\,g_{(m)}=(g_k)_{k=1}^m\in {\mathbb R}^m$
(\cite[Ch.11, \S 6, formulae (7)]{BecBel61})  we have
$$
\Gamma(f_{(m)},g_{(m)})=\sum_{k<n\leq m}M^2(k,n),
$$
therefore, the inequality (\ref{G(fg):G(g)<C.1}) will have the following form
\begin{equation}\label{G(fg):G(g)<C.2}
\frac{\Gamma(f_{(m)},g_{(m)})}{\Gamma(g_{(m)})}=\frac{\sum_{k<n\leq m}M^2(k,n)}{\sum_{k=1}^mg_k^2}\leq C,\quad
m\in {\mathbb N}.
\end{equation}
For $t\in {\mathbb R}$ and $f_{(n)},g_{(n)}\in {\mathbb R}^n$ introduce the function
$$
F_n(t)=\Vert f_{(n)}-tg_{(n)}\Vert^2=(f_{(n)},f_{(n)})-2t(f_{(n)},g_{(n)})+t^2(g_{(n)},g_{(n)}).
$$
The minimum of the function $F_n(t)$ is reached at $t_0^{(n)}=\frac{(f_{(n)},g_{(n)})}{(g_{(n)},g_{(n)})}$ therefore,  we have
\begin{equation*}
F_n(t)=(g_{(n)},g_{(n)})(t-t_0^{(n)})^2+\frac{\Gamma(f_{(n)},g_{(n)})}{\Gamma(g_{(n)})},\quad F_n(t_0^{(n)})=\frac{\Gamma(f_{(n)},g_{(n)})}
{\Gamma(g_{(n)})},
\end{equation*}
hence,
\begin{equation}
\label{F(t_0)-F_{t_n}}
F_n(t_0)-F_n(t_0^{(n)})=(g_{(n)},g_{(n)})(t_0-t_0^{(n)})^2.
\end{equation}
Since $F_n(t_0^{(n)})=\frac{\Gamma(f_{(n)},g_{(n)})}{\Gamma(g_{(n)})}$  is bounded by assumption and
$$
\lim_{n\to\infty}F_n(t)=\lim_{n\to\infty}\Vert f_{(n)}-tg_{(n)}\Vert^2=\infty \quad\text{for all}\quad t\in {\mathbb R},
$$
by the condition (\ref{norm=infty}), we conclude that $\lim_{n\to\infty}\big(F_n(t_0)-F_n(t_0^n)\big)=\infty$.

We show that
condition (\ref{G(fg):G(g)<C.1}) implies that $F_n(t_0)-F_n(t_0^{(n)})$ is bounded.
This contradiction will prove the lemma. Indeed, we have
$$
t_0^{(n+1)}-t_0^{(n)}=\frac{(f_{(n+1)},g_{(n+1)})}{(g_{(n+1)},g_{(n+1)})}-
\frac{(f_{(n)},g_{(n)})}{(g_{(n)},g_{(n)})}=-\frac{\sum_{k=1}^nM(k,n+1)g_kg_{n+1}}
{(g_{(n)},g_{(n)})(g_{(n+1)},g_{(n+1)})}
$$
and
$$
t_0^{(n+m)}-t_0^{(n)}=\frac{(f_{(n+m)},g_{(n+m)})}{(g_{(n+m)},g_{(n+m)})}-
\frac{(f_{(n)},g_{(n)})}{(g_{(n)},g_{(n)})}
$$
\begin{equation}
\label{t(n+m)-t(n)}
=-\frac{\sum_{k=1}^n\sum_{r=n+1}^{n+m}M(k,r)g_kg_r}
{(g_{(n)},g_{(n)})(g_{(n+m)},g_{(n+m)})}=
-\frac{(M_{n,m}g^{n,m},g_{(n)} )}{(g_{(n)},g_{(n)})(g_{(n+m)},g_{(n+m)})}
\end{equation}
where  the vector  $g^{n,m}\in {\mathbb R}^m$ and the rectangular matrix $M_{n,m}\in {\rm Mat}({\mathbb R},n\times m)$ are defined
as follows:
$$g^{n,m}=(g_k)_{k=n+1}^{n+m}\quad\text{ and}\quad M_{n,m}=(M(k,r))_{k,r}\quad 1\leq k\leq n,\quad n+1\leq r\leq n+m.
$$
We observe  that $\lim_{n}t_0^{(n)}=\lim_nt_n=t_0$. Indeed, if $n\to\infty$ by (\ref{alpha(n)-to-0}) we have
$$
t_0^{(n)}=\frac{(f_{(n)},g_{(n)})}{(g_{(n)},g_{(n)})}=\frac{\Vert f_{(n)}\Vert\Vert g_{(n)}\Vert\cos\alpha_n}{\Vert g_{(n)}\Vert^2}=t_n\cos\alpha_n\to t_0.
$$
Finally, for all $n,m\in {\mathbb N}$ we get by (\ref{F(t_0)-F_{t_n}}),  (\ref{t(n+m)-t(n)}) and the Schwartz inequality
$$
F_n(t_0^{(n+m)})-F_n(t_0^{(n)})=(g_{(n)},g_{(n)})(t_0^{(n+m)}-t_0^{(n)})^2=
$$
$$
(g_{(n)},g_{(n)})\left[
\frac{(M_{n,m}g^{n,m},g_{(n)} )}{(g_{(n)},g_{(n)})(g_{(n+m)},g_{(n+m)})}
\right]^2\leq
$$
$$
 \frac{\Vert g_{(n)}\Vert^2 \Vert M_{n,m}g^{n,m} \Vert^2 \Vert g_{(n)}\Vert^2}
{ \Vert g_n{(n)}\Vert^4  \Vert g_{(n+m)}\Vert^4}\leq
\frac{ \Vert M_{n,m}\Vert^2_{\sigma_2}\Vert g^{n,m}\Vert^2}
{\Vert g_{(n+m)}\Vert^4}\leq
\frac{\Vert M_{n+m}\Vert^2_{\sigma_2}}{\Vert g_{(n+m)}\Vert^2}\leq C,
$$
where $M_m:=(M(k,r))_{k<r\leq m}$ and
$$
\Vert M_{n,m}\Vert^2_{\sigma_2}=\sum_{k=1}^n\sum_{r=n+1}^{n+m}M^2(k,r),\quad
\Vert M_{m}\Vert^2_{\sigma_2}=\sum_{k<r\leq m}M^2(k,r)=\Gamma(f_{(m)},g_{(m)}).
$$
Fix $\varepsilon>0$.  Since $\lim_m t_0^{(m)}=t_0$ and the functions $F_n(t)$ are continuous
we conclude that there exists $m_n\geq n$ such that
$
F_n(t_0^{(m)})>F_n(t_0)-\varepsilon, \forall m\geq m_n,
$
in particular, $F_n(t_0^{(m_n)})>F_n(t_0)-\varepsilon$. Since $\lim_nF_n(t_0)=\infty$ we conclude that
$$
\lim_nF_n(t_0^{(m_n)})\geq \lim_n(F_n(t_0)-\varepsilon)=\infty
$$
that contradicts the condition $F_n(t_0^{(n+m)})-F_n(t_0^{(n)})\leq C$ for all $m,n\in {\mathbb N}$.
\qed\end{pf}
\begin{lem}
\label{approx-(x,D)} If $\mu^{L_t}\perp\mu$ for all $t\in
GL(2,{\mathbb R})\setminus\{e\}$, we can approximate one of  the
following  pair of operators:
$(x_{1n},x_{2n}),\,\,(x_{1n},D_{2n}),\,\,(D_{1n},x_{2n}),\,\,$ or
$(D_{1n},D_{2n})$.
\end{lem}
\begin{pf}
{\small
For the convenience of the readers we collect the important formulas below:
\begin{equation}
\label{(xx)1}
\Sigma_{1,m}(x,x)=\frac{\Gamma(f_m^1)+\Gamma(f_m^1,g_m^1)}{\Gamma(g_m^1)+1}=\frac{\sum_{k=-m}^m \frac{b_{1k}^2}{b_{1k}^2+2b_{1k}b_{2k}}
+\Gamma(f_m^1,g_m^1)}
{\sum_{k=-m}^m\frac{b_{2k}^2}{b_{1k}^2+2b_{1k}b_{2k}}+1},
\end{equation}

\begin{equation}
\label{(xx)2}
\Sigma_{2,m}(x,x)=\frac{\Gamma(f_m^2)+\Gamma(f_m^2,g_m^2)}{\Gamma(g_m^2)+1}=\frac{\sum_{k=-m}^m\frac{b_{2k}^2}{b_{2k}^2+2b_{1k}b_{2k}}
+\Gamma(f_m^2,g_m^2)}
{\sum_{k=-m}^m\frac{b_{1k}^2}{b_{2k}^2+2b_{1k}b_{2k}}+1},
\end{equation}
\begin{equation}
\label{(D)1}
\Sigma_{1,m}(D)=\frac{\Gamma(f_m)+\Gamma(f_m,g_m)}{\Gamma(g_m)+1}=\frac{\sum_{k=-m}^m\frac{a_{1k}^2}
{\frac{1}{2b_{1k}}+\frac{1}{2b_{2k}}}+\Gamma(f_m,g_m)}{\sum_{k=-m}^m\frac{a_{2k}^2}
{\frac{1}{2b_{1k}}+\frac{1}{2b_{2k}}}+1},
\end{equation}
\begin{equation}
\label{(D)2}
\Sigma_{2,m}(D)=\frac{\Gamma(g_m)+\Gamma(g_m,f_m)}{\Gamma(f_m)+1}=\frac{\sum_{k=-m}^m\frac{a_{2k}^2}
{\frac{1}{2b_{1k}}+\frac{1}{2b_{2k}}}+\Gamma(g_m,f_m)}{\sum_{k=-m}^m\frac{a_{1k}^2}
{\frac{1}{2b_{1k}}+\frac{1}{2b_{2k}}}+1},
\end{equation}

\begin{equation}
\label{f^1,g^1}
f_m^1=\Big(\frac{b_{1k}}{\sqrt{b_{1k}^2+2b_{1k}b_{2k}}}\Big)_{k=-m}^m,\quad g_m^1=\Big(\frac{b_{2k}}
{\sqrt{b_{1k}^2+2b_{1k}b_{2k}}}\Big)_{k=-m}^m,
\end{equation}
\begin{equation}
\label{f^2,g^2}
f_m^2=\Big(\frac{b_{2k}}{\sqrt{b_{2k}^2+2b_{1k}b_{2k}}}\Big)_{k=-m}^m,\quad g_m^2=\Big(\frac{b_{1k}}{\sqrt{b_{2k}^2+2b_{1k}b_{2k}}}\Big)_
{k=-m}^m,
\end{equation}
\begin{equation}
\label{f_m,g_m=}
f_m=\Big(a_{1k}\Big(\frac{1}{2b_{1k}}+\frac{1}{2b_{2k}}\Big)^{-1/2}\Big)_{k=-m}^m,\quad
g_m=\Big(a_{2k}\Big(\frac{1}{2b_{1k}}+\frac{1}{2b_{2k}}\Big)^{-1/2}\Big)_{k=-m}^m.
\end{equation}

}

To estimate $\Sigma_{1,m}(x,x)$ and $\Sigma_{2,m}(x,x)$ consider three possibilities:
\begin{equation}
(1)\,\, \Sigma^{12}:=\!\sum_{k\in{\mathbb Z}}\frac{b_{1k}}{b_{2k}}\!<\!\infty,\,\,\,
(2)\,\, \Sigma^{21}:=\!\sum_{k\in{\mathbb Z}}\frac{b_{2k}}{b_{1k}}\!<\!\infty,\,\,\,
(3)\,\, \sum_{k\in{\mathbb Z}}\frac{b_{1k}}{b_{2k}}\!=\!\sum_{k\in{\mathbb Z}}\frac{b_{2k}}{b_{1k}}\!=\!\infty.
\end{equation}
 We present the results in the  table I.
 \vskip 0.2cm
\begin{tabular}{|p{1.3cm}|p{2.4cm}|p{2.5cm}|p{1.3cm}|p{1.3cm}|p{1.6cm}|} \hline
\text{table}\,\,I&(1)&(2)&(3a)&(3b)&(3c)\\   \hline
$\Sigma^{12}$&$<\infty$&&$\infty$&$\infty$&$\infty$\\   \hline
$\Sigma^{21}$&&$<\infty$&$\infty$&$\infty$&$\infty$\\   \hline
%
%
$\Vert g^1\Vert$&&&$<\infty$&&\\
  \hline
$\Vert g^2\Vert$&&&&$<\infty$&\\
  \hline
 Lemma &\ref{x2x2},\,\,\ref{d1},     &\ref{x1x1},\,\,\ref{d1},    &\ref{x1x1},&\ref{x2x2},& \ref{x1x1},\,\,\ref{x2x2},\\
       &\ref{d2},\,\,\ref{l.min=proj}&\ref{d2}\,\,\ref{l.min=proj}&   &   & \ref{l.min=proj},\,\,\ref{l.park}  \\
  \hline
&$x_{2n},\,\,D_{1n},D_{2n}$&$x_{1n},\,\,D_{1n},\,\,D_{2n}$&$x_{1n},\,\,x_{2n}$&$x_{1n},\,\,x_{2n}$&$x_{1n},\,\,x_{2n}$\\   \hline
\end{tabular}
 \vskip 0.2cm

{\bf Case (1)}. If $\Sigma^{12}<\infty$ then $\Sigma^{21}=\infty$ and we have $\lim_{m\to\infty}\Sigma_{2,m}(x,x)=\infty$
by Remark~\ref{x1x1>}. Hence, $x_{2n}x_{2t}\,\,\eta\,\,{\mathfrak A}$, by Lemma~\ref{x2x2} and
$x_{2n}\,\,\eta\,\,{\mathfrak A}$, by Lemma~\ref{5.3}. We can approximate $D_{1n}$ and $D_{2n}$ by
Lemmas~\ref{d1}, \ref{d2} and Lemma~\ref{l.min=proj}:
$$
D_{1n}\,\,\eta\,\,{\mathfrak A}\quad\text{if}\,\,\frac{\Gamma(f_m)+\Gamma(f_m,g_m)}{\Gamma(g_m)+1}\to \infty,\,\,
D_{2n}\,\,\eta\,\,{\mathfrak A}\quad\text{if}\,\,\frac{\Gamma(g_m)+\Gamma(g_m,f_m)}{\Gamma(f_m)+1}\to \infty,
$$
where $f_m$ and $g_m$ are defined by (\ref{f_m,g_m=}). Set
\begin{equation}
\label{f,g=}
f=\Big(a_{1k}\Big(\frac{1}{2b_{1k}}+\frac{1}{2b_{2k}}\Big)^{-1/2}\Big)_{k\in{\mathbb Z}},\quad
g=\Big(a_{2k}\Big(\frac{1}{2b_{1k}}+\frac{1}{2b_{2k}}\Big)^{-1/2}\Big)_{k\in{\mathbb Z}}.
\end{equation}
Since $\sum_{k\in{\mathbb Z}}\frac{b_{1k}}{b_{2k}}<\infty$, we
conclude that
\begin{equation}
\label{norm=infty1}
\Vert f\Vert^2=\Vert g\Vert^2=\Vert f-sg\Vert^2=\infty.
\end{equation}
Indeed, we have
$$
\Vert f\Vert^2=\sum_{k\in{\mathbb Z}}\frac{a_{1k}^2}{\frac{1}{2b_{1k}}+\frac{1}{2b_{2k}}}=
\sum_{k\in{\mathbb Z}}\frac{b_{1k}a_{1k}^2}{\frac{1}{2}+\frac{b_{1k}}{2b_{2k}}}\sim 2\sum_{k\in{\mathbb Z}}b_{1k}a_{1k}^2=
S^L_{11}(\mu)=\infty,
$$
$$
\Vert g\Vert^2=
\sum_{k\in{\mathbb Z}}\frac{b_{1k}a_{2k}^2}{\frac{1}{2}+\frac{b_{1k}}{2b_{2k}}}\sim
\sum_{k\in{\mathbb Z}}b_{1k}a_{2k}^2\sim
 \sum_{k\in{\mathbb Z}}\frac{b_{1k}}{2}\Big(\frac{1}{2b_{2k}}+a_{2k}^2\Big)=S^L_{12}(\mu)=\infty,
$$
$$
\Vert f\!-\!sg\Vert^2\!=\!\sum_{k\in{\mathbb Z}}\frac{b_{1k}(a_{1k}-s a_{2k})^2}{\frac{1}{2}+\frac{b_{1k}}{2b_{2k}}}\!\sim\!
\sum_{k\in{\mathbb Z}}b_{1k}(a_{1k}\!-\!s a_{2k})^2\!=\!\frac{1}{4}\sum_{k\in{\mathbb Z}}b_{1k}(-2a_{1k}\!+\!2sa_{2k})^2
$$
$$
\sim
\frac{1}{2}\Big(\frac{(2s)^2}{4}\sum_{k\in{\mathbb
Z}}\frac{b_{1k}}{b_{2k}}+ \sum_{k\in{\mathbb
Z}}\frac{b_{1k}}{2}(-2a_{1k}+2sa_{2k})^2\Big)
= \frac{1}{2}S^{L,-}_{12}(\mu,t)=\infty,
$$
for $t=2s$ (see (\ref{S^(L,-)_(kn)})).
Therefore, by Lemma~\ref{l.min=proj} we conclude (see(\ref{final.})) that
$
\lim_{n\to\infty}\frac{\Gamma(f_{(n)},g_{(n)})}{\Gamma(g_{(n)})}=\infty\quad\text{and}\quad
\lim_{n\to\infty}\frac{\Gamma(f_{(n)},g_{(n)})}{\Gamma(f_{(n)})}=\infty,
$
so $D_{1n},\,\, D_{2n}\,\,\eta\,\,{\mathfrak A} $ by Lemmas~\ref{d1} and \ref{d2}.
Finally, $x_{2n}\,\,D_{1n}\,\,D_{2n}\,\,\eta\,\,{\mathfrak A}$. Now we get
$A_{kn}-x_{2k}D_{2n}=x_{1k}D_{1n},\,\,k,n\in {\mathbb Z}$ and the proof is complete since
we are in the case $m=1$.

{\bf Case (2)}. If $\Sigma^{21}<\infty$ then $\Sigma^{12}=\infty$ and  we have
$\lim_{m\to \infty}\Sigma_{1,m}(x,x)=\infty$,   by Remark~\ref{x1x1>}. Hence,
$x_{1n}x_{1t}\,\eta\,{\mathfrak A}$, by Lemma~\ref{x1x1} and
$ x_{1n}\,\eta\,{\mathfrak A}$, by Lemma~\ref{5.3}. As in the
previous case, the  condition $\sum_{k\in{\mathbb
Z}}\frac{b_{2k}}{b_{1k}}<\infty$ implies
$$
\Vert f\Vert^2\sim S^L_{21}(\mu)=\infty,\quad \Vert g\Vert^2\sim S^L_{22}(\mu)=\infty,\quad
\Vert f-sg\Vert^2\sim S^{L,-}_{21}(\mu,t)=\infty,
$$
for $t=\frac{2}{s} $. Exactly, as in the case (1), we can approximate $D_{1n}$ and $D_{2n}$.
Finally, $x_{1n}\,\,D_{1n}\,\,D_{2n}\,\,\eta\,\,{\mathfrak A}$.
Further, $A_{kn}-x_{1k}D_{1n}=x_{2k}D_{2n},\,\,k,n\in {\mathbb Z}$ and the proof is complete.

{\bf Case (3)}. Let $\sum_{k\in{\mathbb Z}}\frac{b_{1k}}{b_{2k}}\!=\!\sum_{k\in{\mathbb Z}}\frac{b_{2k}}{b_{1k}}\!=\!\infty.$
Set $c_n=\frac{b_{2n}}{b_{1n}},\,\,n\in{\mathbb Z}$.
 The vectors $f_m^1,\,g_m^1,\,f_m^2,\,g_m^2$ are defined as follows (see (\ref{f^1,g^1}) and (\ref{f^2,g^2})):
\begin{eqnarray}
\label{f^1,g^1,1}
f_m^1=\Big(\frac{1}{\sqrt{1+2c_n}}\Big)_{-m}^m,\quad&
g_m^1=\Big(\frac{c_n}{\sqrt{1+2c_n}}\Big)_{-m}^m,\\
\label{f^2,g^2,2}
f_m^2=\Big(\sqrt{\frac{c_n}{c_n+2}}\Big)_{-m}^m,\quad&
g_m^2=\Big(\frac{1}{\sqrt{c_n^2+2c_n}}\Big)_{-m}^m.
\end{eqnarray}
We show that
\begin{equation}
\label{f1,f2,g1+g2=inf}
\Vert f^1\Vert^2=\Vert f^2\Vert^2=\Vert g^1\Vert^2+\Vert g^2\Vert^2=\infty.
\end{equation}
Indeed, we have
$$
\Vert f^1\Vert^2=\sum_{n\in {\mathbb Z}}
(1+2c_n)^{-1}
\sim \sum_{n\in {\mathbb Z}}c_n^{-1}=\Sigma^{12}=\infty,
$$
$$
\Vert f^2\Vert^2=\sum_{n\in {\mathbb Z}}
c_n(c_n+2)^{-1}
\sim \sum_{n\in {\mathbb Z}}c_n=\Sigma^{21}=\infty.
$$
Let us suppose that $\Vert g^1\Vert^2+\Vert g^2\Vert^2<\infty$ then
$$
\infty>\Vert g^1\Vert^2+\Vert g^2\Vert^2=\sum_{n\in {\mathbb Z}}\Big(\frac{c_n^2}{1+2c_n}+\frac{1}{c_n^2+2c_n}\Big)>
\sum_{n\in {\mathbb Z}}\frac{1+c_n^2}{(1+c_n)^2},
$$
hence, $\sum_{n\in {\mathbb Z}}\frac{1}{(1+c_n)^2}<\infty$ and $\sum_{n\in {\mathbb Z}}\frac{c_n^2}{(1+c_n)^2}<\infty$
therefore,
$$
\infty>\sum_{n\in {\mathbb Z}}\frac{(1+c_n)^2}{(1+c_n)^2}=\sum_{n\in {\mathbb Z}}1=\infty.
$$
This contradiction proves that $\Vert g^1\Vert^2+\Vert g^2\Vert^2=\infty$. We shall come back to the case $I(3)$ later.
{\it We show that  in the case $A$  (see (\ref{A,B})) we can approximate $x_{1n}$ and $x_{2n}$}.

Now we study the possibility of the approximation of $D_{1n}$ and $D_{2n}$ by Lemmas ~\ref{d1}, \ref{d2} and
\ref{l.min=proj}. Recall the notations:
\begin{equation}
\label{I}
\Vert f_m\Vert^2\!=\!
\sum_{k=-m}^ma_{1k}^2\Big(\frac{1}{2b_{1k}}+\frac{1}{2b_{2k}}\Big)^{-1}\!,\quad
\Vert g_m\Vert^2\!=\!\sum_{k=-m}^ma_{2k}^2\Big(\frac{1}{2b_{1k}}+\frac{1}{2b_{2k}}\Big)^{-1}\!.
\end{equation}

All the different cases are presented in the following tables:
%
 \vskip 0.1cm
\begin{tabular}{|p{1.2cm}|p{1.4cm}|p{1.4cm}|p{1.4cm}|p{1.4cm}|p{2.1cm}|p{0.9cm}|
} \hline
\text{table}\,\,II &(1)&(2)&(3a)&(3b)&(3c)&(4)
\\   \hline
$\Vert f \Vert^2$&$\infty$&$<\infty$&$\infty$&$\infty$&$\infty$&$<\infty$
\\   \hline
$\Vert g \Vert^2$&$<\infty$&$\infty$&$\infty$&$\infty$&$\infty$&$<\infty$
\\   \hline
$
\frac{\Vert f_m\Vert^2}{\Vert g_m\Vert^2}
$&&&$\to\infty$&$\to 0$&
$C_1\!\!\leq\!\!\frac{\Vert f_m\Vert^2}{\Vert g_m\Vert^2}\!\!\leq\!\!C_2$&
\\   \hline
Lemma &\ref{d1} &\ref{d2}&\ref{d1}   &\ref{d2} &\ref{d1} ,\,\,\ref{d2}   &
\\
&\ref{x1,dA} &\ref{x2,dA}&\ref{x1,dA}&\ref{x2,dA}&\ref{the-last},\,\,\ref{l.min=proj}  &   \\
\hline
&$D_{1n},\,\,x_{1n}$&$D_{2n},\,\,x_{2n}$&$D_{1n},\,\,x_{1n}$&$D_{2n},\,\,x_{2n}$&$D_{1n},\,\,D_{2n}$&
\\   \hline
\end{tabular}
 \vskip 0.3cm
\begin{rem}
\label{r.II.1} We show that if $\Vert g \Vert^2<\infty$ and
$S^L_{12}(\mu)=\infty$ then $\sum_n\frac{b_{1n}}{b_{2n}}=\infty$.
Indeed, let us suppose that $\sum_n\frac{b_{1n}}{b_{2n}}<\infty$,
then
\begin{equation}
\label{Ineq.1}
\Vert g \Vert^2=\sum_{n\in{\mathbb Z}}\frac{a_{2n}^2}{\frac{1}{2b_{1n}}+\frac{1}{2b_{2n}}}
\sim
\sum_{n\in{\mathbb Z}}b_{1n}a_{2n}^2
\sim
\sum_{n\in{\mathbb Z}}\frac{b_{1n}}{2}\Big(\frac{1}{2b_{2n}}+a_{2n}^2\Big)=S_{12}^L(\mu)=\infty.
\end{equation}
\end{rem}
We explain the tables ${\rm II}$ in details. {\it The first two case (1) and (2) are independent of the case ${\rm I(3)}$}.

{\bf (1)} If $\Vert g \Vert^2<\infty$ and $\Vert f \Vert^2=\infty$, we
have  $D_{1k}\,\,\eta\,\,{\mathfrak A}$ by Lemma \ref{d1}. The
condition $\Vert g \Vert^2<\infty$ implies $\sum_{k\in{\mathbb
Z}}\frac{b_{1k}}{b_{2k}}=\infty$, by  Remark~\ref{r.II.1}
therefore, $x_{1k}\,\,\eta\,\,{\mathfrak A}$,  by Lemma
\ref{x1,dA}. Further, $A_{kn}-x_{1k}D_{1n}=x_{2k}D_{2n},\,\,k,n\in
{\mathbb Z}$ and the proof is complete since we are reduced to the
case $m=1$.

{\bf (2)} If $\Vert g \Vert^2=\infty$ and $\Vert f \Vert^2<\infty$, we
have  $D_{2k}\,\,\eta\,\,{\mathfrak A}$ by Lemma~\ref{d2}.
 By remark similar to the Remark~\ref{r.II.1}, we conclude that $\sum_{k\in{\mathbb Z}}\frac{b_{2k}}{b_{1k}}=\infty$ therefore,
 $x_{2k}\,\,\eta\,\,{\mathfrak A}$ by Lemma~\ref{x2,dA} and $A_{kn}-x_{2k}D_{2n}=x_{1k}D_{1n},\,\,k,n\in {\mathbb Z}$, case $m=1$.

{\bf (3)} {\it Consider now the case} ${\rm I(3)}$. Let both series be divergent: $\Vert g \Vert^2\!=\!\infty$ and
$\Vert f \Vert^2\!=\!\infty$. {\it We show that in the case $(B)$   (see (\ref{A,B})) holds $\Vert f+sg\Vert^2\!=\!\infty$ for all
$s\in {\mathbb R}$, by Lemma~\ref{l.park} therefore, by Lemma~\ref{l.min=proj}, we can approximate $D_{1n}$ and $D_{2n}$}.
To be more precise  consider three possibilities:

{\bf (3a)} let $\frac{\Vert f_m\Vert^2}{\Vert g_m\Vert^2}\to\infty$, then
$D_{1k}\,\,\eta\,\,{\mathfrak A}$. Since
$\sum_n\frac{b_{1n}}{b_{2n}}=\infty$, we have
$x_{1n}\,\,\eta\,\,{\mathfrak A}$ by  Lemma \ref{x1,dA} and
finally, $x_{1n},\,D_{1n}\,\,\eta\,\,{\mathfrak A},\,\,n\in
{\mathbb Z}$. We are reduced to the case $m=1$.

{\bf (3b)}  Let $\frac{\Vert f_m\Vert^2}{\Vert g_m\Vert^2}\to 0$, then
$D_{2k}\,\,\eta\,\,{\mathfrak A}$. Since
$\sum_n\frac{b_{2n}}{b_{1n}}=\infty$, we get
$x_{2n}\,\,\eta\,\,{\mathfrak A}$, by  Lemma \ref{x2,dA} and
finally, $x_{2n},\,D_{2n}\,\,\eta\,\,{\mathfrak A},\,\,n\in
{\mathbb Z}$. We are reduced to the case $m=1$.

{\bf (3c)} The case when  $\Vert f \Vert^2=\Vert g \Vert^2=\infty$ and $C_1\!\!\leq\!\!
\frac{\Vert f_m\Vert^2}{\Vert g_m\Vert^2}\!\!\leq\!\!C_2$ .

{\bf (4)} The case when $\Vert f \Vert^2+\Vert g \Vert^2<\infty.$

{\it To complete the proof} of the lemma it remains to consider ${\rm I(3)}$, i.e.,  the last  case ${\rm (3)}$
in the table ${\rm I}$ and the last two  cases in the table ${\rm
II}$, i.e., ${\rm II(3c)}$ and  ${\rm II(4)}$, where:
\begin{equation}
\label{I(3c3)}
{\rm I(3)}\quad
\sum_{k\in{\mathbb Z}}\frac{b_{1k}}{b_{2k}}=\sum_{k\in{\mathbb Z}}\frac{b_{2k}}{b_{1k}}=\infty,
\end{equation}
\begin{equation}
{\rm II(c3)}\quad \sum_{k\in{\mathbb Z}}a_{1k}^2\Big(\frac{1}{2b_{1k}}+\frac{1}{2b_{2k}}\Big)^{-1}=
\sum_{k\in{\mathbb Z}}a_{2k}^2\Big(\frac{1}{2b_{1k}}+\frac{1}{2b_{2k}}\Big)^{-1}=\infty,
%
\end{equation}
\begin{equation}
\label{I+II(d)}
{\rm II(4)}\quad
%
\sum_{k\in{\mathbb Z}}\Big(a_{1k}^2+a_{2k}^2\Big)\Big(\frac{1}{2b_{1k}}+\frac{1}{2b_{2k}}\Big)^{-1}<\infty.
%
\end{equation}
Come back to the condition $\mu^{L_t}\perp\mu$. By Remark~\ref{perp2-1} we have
$$
\mu^{L_{\tau_-(\phi,s)}}\perp\mu,\,\,\,\phi\in[0,2\pi),\,\,s>0\Leftrightarrow
\Sigma_1(s)+\Sigma_2(C_1,C_2)\!=\!\infty,\,\,s>0,
$$
for $(C_1,C_2)\in{\mathbb R}^2\setminus\{0\}$.
Recall that (see (\ref{sigma2(s)}))
$$
\Sigma_1(s)=\sum_{n\in {\mathbb Z}}
\Big(s^2\sqrt{\frac{b_{1n}}{b_{2n}}}-s^{-2}\!\sqrt{\frac{b_{2n}}{b_{1n}}}\Big)^2,
$$
$$
\Sigma_2(C_1,C_2)=\sum_{n\in {\mathbb Z}}(C^2_1b_{1n}+C^2_2b_{2n})(C_1a_{1n}+C_2a_{2n})^2.
$$
The condition $\Sigma_1(s)+\Sigma_2(C_1,C_2)\!=\!\infty$ splits into two cases:
\begin{equation}
\label{A,B}
\begin{array}{cccc}
(A) & \Sigma_1(s)=\infty,&&\\
(B) &\Sigma_1(s)<\infty&\text{ but}& \Sigma_2(C_1,C_2)=\infty.
\end{array}
\end{equation}
${\bf (A)\&I(3)}$.  In this case
independently of the conditions ${\rm II(3c)}$ and ${\rm II(4)}$ we can approximate
$x_{1n}$ and $x_{2n}$ by Lemma \ref{x1x1} and \ref{x2x2}.\\
${\bf (B)\&II(3c)}$ In this case
we can  approximate $D_{1n}$ and $D_{2n}$ by Lemmas~\ref{d1} and
\ref{d2} respectively. More precisely, to use Lemma~\ref{l.min=proj} we show that conditions
(\ref{norm=infty}) are satisfied for two vectors $f$ and $g$ defined by (\ref{f_m,g_m=}) (see Lemma~\ref{the-last}).\\
${\bf (B)\&II(4)}$  This  case
(see (\ref{I+II(d)})) can not be realized if $\Sigma_2(C_1,C_2)=\infty$.

{\bf Case} ${\rm (A)\&I(3)}$. Using Lemma~\ref{l.min=proj} we  conclude that
 \begin{equation}
\label{xx3}
\Gamma(f_m^1,g_m^1)(\Gamma(g_m^1))^{-1}\to\infty\quad\text{and}\quad \Gamma(f_m^2,g_m^2)(\Gamma(g_m^2))^{-1}\to\infty.
\end{equation}
To use Lemma~\ref{l.min=proj}, it is sufficient to show that in the case $(A)$ relations (\ref{norm=infty}) hold for
$f^1,g^1$ and $f^2,g^2$, i.e., for all
$s\in {\mathbb R}\setminus\{0\}$ we have (see Lemma~\ref{l.park})
\begin{equation}
\label{norm=infty(1,2)}
\Vert f^1\Vert^2\!=\!\Vert g^1\Vert^2\!=\!\Vert f^1+s g^1\Vert^2=\infty,\,\,\,
\Vert f^2\Vert^2\!=\!\Vert g^2\Vert^2\!=\!\Vert f^2+s g^2\Vert^2=\infty.
\end{equation}
Consider three possibilities in the case ${\rm I(3)}$:\\
%

(3a) If $\Vert g^1\Vert<\infty$ then $\Vert g^2\Vert=\infty$ therefore, we have $\Vert f_m^1\Vert /\Vert g_m^1\Vert\to\infty$ so,
$x_{1n}\,\,\eta\,\,{\mathfrak A}$ by Lemma~\ref{l.min=proj0} (a). In the case $(A)$ by Lemma~\ref{l.park} holds
$\Vert f^2\Vert^2=\Vert g^2\Vert^2=\Vert f^2+s g^2\Vert^2=\infty$ therefore, $x_{2n}\,\,\eta\,\,{\mathfrak A}$ by Lemma\ref{l.min=proj}.

(3a) If $\Vert g^2\Vert<\infty$ then $\Vert g^1\Vert=\infty$ therefore, we have $\Vert f_m^2\Vert /\Vert g_m^2\Vert\to\infty$ so,
$x_{2n}\,\,\eta\,\,{\mathfrak A}$ by Lemma\ref{l.min=proj0} (a). In the case $(A)$ by Lemma~\ref{l.park} holds
$\Vert f^1\Vert=\Vert g^1\Vert=\Vert f^1+s g^1\Vert=\infty$ therefore, $x_{1n}\,\,\eta\,\,{\mathfrak A}$ by Lemma~\ref{l.min=proj}.

(3c) If $\Vert g^1\Vert=\Vert g^2\Vert=\infty$ then by Lemma~\ref{l.park} all relations (\ref{norm=infty(1,2)}) hold in the case $(A)$
therefore, $x_{1n},\,\,x_{2n}\,\,\eta\,\,{\mathfrak A}$.

To prove (\ref{norm=infty(1,2)}) we need the following auxiliary lemma.
\begin{lem}
\label{l.2-equiv-perp}
 The  following two conditions are equivalent:
\begin{equation}
\label{2-equiv-perp1}
(i)\quad\Sigma_1(s)=\sum_{n\in {\mathbb Z}}
\Big(s^2\sqrt{\frac{b_{1n}}{b_{2n}}}-s^{-2}\sqrt{\frac{b_{2n}}{b_{1n}}}\Big)^2=\infty,
\end{equation}
\begin{equation}
\label{2-equiv-perp2}
(ii)\quad\Sigma_2(s)=\sum_{n\in {\mathbb Z}}
\Big(s^4\frac{b_{1n}}{b_{2n}}-1\Big)^2+\Big(s^{-4}\frac{b_{2n}}{b_{1n}}-1\Big)^2=\infty.
\end{equation}
\end{lem}

\begin{pf}
We show that $(i)\Rightarrow (ii)$. Indeed, we have
$$
(a^2-1)^2+(a^{-2}-1)^2=(a^2-1)^2(1+a^{-4})=(a-a^{-1})^2(a^2+a^{-2}).
$$
Set $a=s^2(b_{1n}/b_{2n})^{1/2}$, then
$$
\Sigma_2(s)=\sum_{n\in {\mathbb Z}}
\Big(s^2\sqrt{\frac{b_{1n}}{b_{2n}}}-s^{-2}\sqrt{\frac{b_{2n}}{b_{1n}}}\Big)^2
\Big(s^4\frac{b_{1n}}{b_{2n}}+s^{-4}\frac{b_{2n}}{b_{1n}}\Big)\geq 2\Sigma_1(s).
$$
We prove that $(ii)\Rightarrow (i)$. Denote by $s^4\frac{b_{1n}}{b_{2n}}=1+a_n$, then we have
$$
\Sigma_1(s)=\sum_{n\in {\mathbb Z}}\Big(\sqrt{1+a_n}-\frac{1}{\sqrt{1+a_n}}\Big)^2=\sum_{n\in {\mathbb Z}}
\Big(\frac{a_n}{\sqrt{1+a_n}}\Big)^2=\sum_{n\in {\mathbb Z}}\frac{a_n^2}{1+a_n},
$$
$$
\Sigma_2(s)\!=\!\sum_{n\in {\mathbb
Z}}\Big(a_n^2+\Big(\frac{1}{1+a_n}\!-\!1\Big)^2\Big)\!=\!\sum_{n\in {\mathbb
Z}}\Big(a^2\!+\!\frac{a_n^2}{(1+a_n)^2}\Big)\!\!
\stackrel{(\ref{sim1})}{\sim}\!\! \sum_{n\in {\mathbb
Z}}a_n^2+\sum_{n\in {\mathbb Z}}\frac{a_n^2}{1+a_n}.
$$
Let $\Sigma_2(s)=\infty$. If $\sum_{n\in {\mathbb
Z}}\frac{a_n^2}{1+a_n}=\infty$, the proof is  complete. Suppose
that $\sum_{n\in {\mathbb Z}}a^2_n=\infty$. We show that in this
case $\Sigma_1(s)=\infty$. It is sufficient to prove that
$$
\sum_{n\in {\mathbb N}}a_n^2=\infty \quad\text{implies}\quad \sum_{n\in {\mathbb N}}a_n^2(1+a_n)^{-1}=\infty.
$$
Consider three cases:\\
(a) If $0<\varepsilon \leq  1+a_n\leq C<\infty$ for all $n\in
{\mathbb N}$, then
$$
C^{-1}\sum_{n\in {\mathbb N}}a_n^2\leq \sum_{n\in {\mathbb N}}a_n^2(1+a_n)^{-1}\leq \varepsilon^{-1}\sum_{n\in {\mathbb N}}a_n^2.
$$
(b) If $\lim_{k\to\infty}(1+a_{n_k})=0$, then
$$
\lim_{k\to\infty} a_{n_k}^2(1+a_{n_k})^{-1}=\infty\quad\text{and}\quad \sum_{n\in {\mathbb N}}a_n^2(1+a_n)^{-1}=\infty.
$$
(c) If $\lim_{k\to\infty}(1+a_{n_k})=+\infty$, then
$$
\sum_{n\in {\mathbb N}}a_n^2(1+a_n)^{-1}>\sum_{k\in {\mathbb N}}a_{n_k}(a_{n_k}^{-1}+1)^{-1}\sim \sum_{k\in {\mathbb N}}a_{n_k}=\infty.
$$
\qed\end{pf}

\begin{lem}
\label{l.park} If $\Sigma_1(s)=\infty$ for any $s>0$, then
$$
\Vert f^1-Cg^1\Vert^2=\infty\quad\text{and}\quad \Vert f^2-Cg^2\Vert^2=\infty,\quad\text{for any}\quad C>0.
$$
\end{lem}
\begin{pf}
Set as before $c_n=\frac{b_{2n}}{b_{1n}},\,\,n\in{\mathbb Z}$. Suppose that $\Sigma_1(s)=\infty$, then
$$
\infty=\Sigma_1(s)=\sum_{n\in{\mathbb Z}}\Big(\frac{s^2}{\sqrt{c_n}}-\frac{\sqrt{c_n}}{s^2}\Big)^2=
\sum_{n\in {\mathbb Z}}\frac{a_n^2}{1+a_n},
$$
where $s^4c_n^{-1}=1+a_n$ or $c_n=\frac{s^4}{1+a_n}$. We show that
$$
\Vert s^4f^1-g^1\Vert^2=\infty
\quad\text{and}\quad
\Vert s^{-4}f^2-Cg^2\Vert^2
=\infty.
$$
Indeed, using (\ref{f^1,g^1,1}) and (\ref{f^2,g^2,2}) we get
$$
\Vert s^4f^1-g^1\Vert^2=\sum_{k\in{\mathbb Z}}\frac{(s^4-c_k)^2}{1+2c_k}=
\sum_{k\in{\mathbb Z}} \Big(\frac{s^4}{c_k}-1\Big)\Big(\frac{1}{c_k^2}+\frac{2}{c_k}\Big)^{-1}=
$$
$$
\sum_{k\in{\mathbb Z}} \frac{a_k^2
}{\big(\frac{1+a_k}{s^4}\big)^2+2\frac{1+a_k}{s^4}}\sim
\sum_{k\in{\mathbb Z}} \frac{a_k^2
}{
(1+a_k)^2+2(1+a_k)
}\!=\!
$$
$$
\sum_{k\in{\mathbb Z}}\frac{a_k^2}{3+4a_k+a_k^2}\sim \sum_{k\in{\mathbb Z}} \frac{a_k^2}{1+a_k}=\infty
$$
and
$$
\Vert s^{-4}f^2-g^2\Vert^2=\sum_{k\in{\mathbb Z}} \frac{(s^{-4}c_k^2-1)^2}{c_k^2+2c_k}=
$$
$$
\sum_{k\in{\mathbb Z}} \Big(\frac{1}{1+a_k}-1\Big)^2
\Big(\Big(\frac{s^4}{1+a_k}\Big)^2+2\frac{s^4}{1+a_k}\Big)^{-1}
=
$$
$$
\sum_{k\in{\mathbb Z}} \frac{a_k^2
}{
s^8+2s^4(1+a_k)
}\sim \sum_{k\in{\mathbb Z}}
 \frac{a_k^2}{1+a_k}=\infty.
$$
\qed\end{pf}
So, in the case ${\rm (A)\&I(3)}$  we can approximate $x_{1n}$ and $x_{2n}$.

{\bf Case} ${\rm (B)\&II(3c)}$.
\begin{lem}
\label{the-last} When $\Sigma_1(s)<\infty$ and
$\Sigma_2(C_1,C_2)=\infty$,   we get
\begin{equation}
\label{C_1f+C_2g=}
\sigma(C_1,C_2):=
\Vert C_1f+C_2g\Vert^2=\sum_{n\in {\mathbb Z}}\frac{(C_1a_{1n}+C_2a_{2n})^2}{\frac{1}{2b_{1n}}+\frac{1}{2b_{2n}}}=\infty,\,\,
(C_1,C_2)\in{\mathbb R}^2\setminus\{0\},
\end{equation}
where $f$ and $g$ are defined by (\ref{f,g=})
$$
f=\Big(a_{1k}\Big(\frac{1}{2b_{1k}}+\frac{1}{2b_{2k}}\Big)^{-1/2}\Big)_{k\in{\mathbb Z}},\quad
g=\Big(a_{2k}\Big(\frac{1}{2b_{1k}}+\frac{1}{2b_{2k}}\Big)^{-1/2}\Big)_{k\in{\mathbb Z}}.
$$
\end{lem}
\begin{pf}
Let $\Sigma_1(s)=\sum_{n\in {\mathbb
Z}}\frac{a_n^2}{1+a_n}<\infty$, where
$s^4\frac{b_{1n}}{b_{2n}}=1+a_n$ or $s^4b_{1n}= (1+a_n)b_{2n}$. We
see that $\lim_n\frac{a_n^2}{1+a_n}=0$ hence,
$\lim_na_n=\lim_n\big(s^4\frac{b_{1n}}{b_{2n}}-1\big)=0$. We have
$$
\sigma(C_1,C_2)=\sum_{n\in {\mathbb Z}}\frac{b_{1n}(C_1a_{1n}+C_2a_{2n})^2}{\frac{1}{2}+\frac{b_{1n}}{2b_{2n}}}
=\sum_{n\in {\mathbb Z}}\frac{b_{1n}(C_1a_{1n}+C_2a_{2n})^2}{\frac{1}{2}+\frac{1}{2}\frac{1+a_n}{s^4}}
$$
$$
\sim \sum_{n\in {\mathbb Z}}C^2_1b_{1n}(C_1a_{1n}+C_2a_{2n})^2,
$$
$$
\sigma(C_1,C_2)=\sum_{n\in {\mathbb Z}}\frac{b_{2n}(C_1a_{1n}+C_2a_{2n})^2}{\frac{b_{2n}}{2b_{1n}}+\frac{1}{2}}
=\sum_{n\in {\mathbb Z}}\frac{b_{2n}(C_1a_{1n}+C_2a_{2n})^2}{\frac{1}{2}\frac{s^4}{1+a_n}+\frac{1}{2}}
$$
$$
\sim \sum_{n\in {\mathbb Z}}C^2_2b_{2n}(C_1a_{1n}+C_2a_{2n})^2,
$$
hence, $\sigma(C_1,C_2)\sim  \sum_{n\in {\mathbb Z}}(C^2_1b_{1n}+C^2_2b_{2n})(C_1a_{1n}+C_2a_{2n})^2=
\Sigma_2(C_1,C_2).$
\qed\end{pf}
Finally, we can approximate $D_{1n}$ and $D_{2n}$ in the case ${\rm (B)\&II(3c)}$.

{\bf Case} ${\rm (B)\&II(4)}$. The last case  ${\rm (B)\&II(4)}$ (see (\ref{I+II(d)})) can not be realized if
$\Sigma_2(C_1,C_2)=\infty$.
Indeed, in this case by Lemma~\ref{the-last} $ \sigma(C_1,C_2)\sim \Sigma_2(C_1,C_2)=\infty$.
This contradicts (\ref{I+II(d)}):
 $$
 \sum_{k\in{\mathbb Z}}\big(a_{1k}^2+a_{2k}^2\big)
\Big(\frac{1}{2b_{1k}}+\frac{1}{2b_{2k}}\Big)^{-1}<\infty.
$$ This completes the proof of Lemma~\ref{approx-(x,D)} for $m=2$.
\qed\end{pf}
The proof of  the irreducibility for $m=2$ follows  from
Remark~\ref{r.irr-idea-G}. Depending on the measure, we can
approximate four different families of commuting operators
$B^\alpha=(B^\alpha_{1n},B^\alpha_{2n})_{n\in{\mathbb Z}}$ for
$\alpha\in \{0,1\}^2$:
$$
B^{(0,0)}\!\!=\!\!(x_{1n},x_{2n})_{n},\,\,B^{(0,1)}\!\!=\!\!(x_{1n},D_{2n})_{n},\,\,
B^{(1,0)}\!\!=\!\!(D_{1n},x_{2n})_{n},\,\,B^{(0,0)}\!\!=\!\!(D_{1n},D_{2n})_{n}.
$$
The von Neumann algebra $L^\infty_\alpha(X_2,\mu^2)$ consists of
all essentially bounded functions $f(B^\alpha)$ in the commuting
family of operators $B^\alpha$ (see, e.g., \cite{Ber86}) as, in
particular, $L^\infty_{(0,0)}(X_2,\mu^2)=L^\infty(X_2,\mu^2)$.
Since the von Neumann algebras $L^\infty_\alpha(X_2,\mu^2)$ are
maximal abelian, the commutant $\big({\mathfrak A}^2\big)'$ of the
von Neumann algebra ${\mathfrak A}^2$ generated by the
representation is contained in  $L^\infty_\alpha(X_2,\mu^2)$.
Hence, the bounded operator $A\in \big({\mathfrak A}^2\big)'$ will
be some function $A=a(B^\alpha)\in L^\infty_\alpha(X_2,\mu^2)$.
The commutation relation $[A,T^{R,\mu,2}_t]=0$ gives us the
following relations:
$a((B^\alpha)^{R_t})=a(B^\alpha)$ for all $t\in {\rm
GL}_0(2\infty,{\mathbb R}).$
Set $B^\alpha_r=(B^\alpha_{rn})_n,\,\,x_r=(x_{rn})_n,\,\,D_r=(x_{rn})_n,\,\,r=1,2,\,\,n\in{\mathbb
Z}$ and set as before, $E_{kn}(t):=I+tE_{kn},\,\,t\in {\mathbb R},\,\,k,n\in
{\mathbb Z},\,\,k\not=n$. Then the action $(B^\alpha)^{R_s}$ is
defined as follows:
%
\begin{eqnarray*}
(B^\alpha_1,B^\alpha_2)^{R_t}=((B^\alpha_1)^{R_t},(B^\alpha_2)^{R_t}),\quad (x_r)^{R_t}=x_rt,\quad (D_r)^{R_t}=D_rt^*,\\
a(\dots, x_{rk},\dots ,x_{rn},\dots)^{R_{E_{kn}(t)}}=a(\dots, x_{rk},\dots ,x_{rn}+tx_{rk},\dots),\\
a(\dots, D_{rk},\dots ,D_{rn},\dots)^{R_{E_{kn}(t)}}=a(\dots, D_{rk}+tD_{rn},\dots ,D_{rn},\dots),\,\,\,t\in {\mathbb R}.
\end{eqnarray*}
In all the cases, by ergodicity of the measure $\mu^2$, we conclude that $a$ is constant.
\subsection{The proof of  Lemmas  \ref{perp1}, \ref{perp2}}
\label{ch5.6}  Lemma \ref{perp1} follows from Lemmas \ref{perp3}-
\ref{ldet*}.
\begin{lem}
\label{perp3} For $t\in
 {\rm GL}(m,{\mathbb R})\backslash\{e\}$ we have $ (\mu_{(b,a)}^m)^{L_t}\perp\mu_{(b,a)}^m
 $if and only if
 \begin{equation}
  \label{7.1}
(\mu_{(b,0)}^m)^{L_t}\perp\mu_{(b,0)}^m\quad\text{or}\quad
\mu_{(b,{L_t}a)}^m\perp\mu_{(b,a)}^m.
 \end{equation}
\end{lem}
Let us define the following measures on the spaces ${\mathbb R}^m$
and $X_m$:
$$
\mu_m^{(B_n,0)}=\otimes_{k=1}^m\mu_{(b_{kn},0)},\quad
\mu_m^{(B_n,a_n)}=\otimes_{k=1}^m\mu_{(b_{kn},a_{kn})},
$$
where $a_n\!=\!(a_{1n},...,a_{mn})\in{\mathbb R}^m$ and
$B_n\!=\!{\rm diag}(b_{1n},...,b_{mn})\in {\rm Mat}(m,{\mathbb R})$. Since
$$
\mu_{(b,a)}^m=\otimes_{n\in{\mathbb Z}}\mu_m^{(B_n,a_n)},\quad
\mu_{(b,0)}^m=\otimes_{n\in{\mathbb Z}}\mu_m^{(B_n,0)},
$$
$$
\left(\mu_{(b,a)}^m\right)^{L_t}=\otimes_{n\in{\mathbb
Z}}\left(\mu_m^{(B_n,a_n)}\right)^{L_t},\quad
\left(\mu_{(b,0)}^m\right)^{L_t}=\otimes_{n\in{\mathbb
Z}}\left(\mu_m^{(B_n,0)}\right)^{L_t},
$$
and
$$
\mu_{(b,{L_t}a)}^m=\otimes_{n\in{\mathbb
Z}}\mu_m^{(B_n,{L_t}a_n)},
$$
by Kakutani criterion \cite{Kak48}, we have two lemmas:

\begin{lem}
\label{l7.2}
 For measures $\mu_{(b,0)}^m,\,m\in{\mathbb N}$ and $t\in
 {\rm GL}(m,{\mathbb R})\backslash\{e\}$ we obtain
$$
(\mu_{(b,0)}^m)^{L_t}\perp\mu_{(b,0)}^m\,\, \Leftrightarrow
\prod_{n\in{\mathbb
Z}}H\left(\left(\mu_m^{(B_n,0)}\right)^{L_t},\mu_m^{(B_n,0)}\right)=0.
$$
\end{lem}
\begin{lem}
\label{l7.3}
 For measures $\mu_{(b,0)}^m,\,m\in{\mathbb N}$ and $t\in
 {\rm GL}(m,{\mathbb R})\backslash\{e\}$ we get
$$
\mu_{(b,{L_t}a)}^m\perp\mu_{(b,a)}^m \Leftrightarrow
\prod_{n\in{\mathbb
Z}}H\left(\mu_m^{(B_n,{L_t}a_n)},\mu_m^{(B_n,a_n)}\right)=0.
$$
\end{lem}
To prove Lemma \ref{perp1} it is sufficient to show, due to Lemma
\ref{perp3}, that
\begin{equation}
\label{Hel1-sl}
H\left(\left(\mu_m^{(B_n,0)}\right)^{L_t},\mu_m^{(B_n,0)}\right)=
\left(\frac{1}{2^m\vert{\rm det}\,\,t\vert} {\rm
det}\left(I+X_n^*(t)X_n(t)\right)\right)^{-1/2},
\end{equation}
to prove the  equivalence
\begin{equation}
\label{Hel2-sl}
 \prod_{n\in{\mathbb
Z}}H\left(\mu_m^{(B_n,{L_t}a_n)},\mu_m^{(B_n,a_n)}\right)=0
\Leftrightarrow \sum_{n\in {\mathbb Z}}\sum_{r=1}^m
b_{rn}\Big(\sum_{s=1}^m(t_{rs}-\delta_{rs})a_{sn}\Big)^2 =\infty,
\end{equation}
and to use the following lemma:
\begin{lem}
\label{ldet*}
 For $X\in {\rm Mat}(m,{\mathbb R})$ we have
\begin{equation}
\label{det*} {\rm det}\left(I+X^*X\right)= 1+\sum_{r=1}^{m}
\sum_{1\leq i_1<i_2<...<i_r\leq m;1\leq j_1<j_2<...<j_r\leq m}
\left(M^{i_1i_2...i_r}_{j_1j_2...j_r}(X)\right)^2.
\end{equation}
\end{lem}
{\it The proof} of  equality (\ref{Hel1-sl}) is based on the exact
formula of the Hellinger integral  (see \cite{Kuo75} for
definition) for two  Gaussian measures $\mu=\mu_m^{(B_n,0)}$ and
$\nu=\mu_m^{(C_n,0)}$ in the space ${\mathbb R}^m$ (see
\cite{Kuo75})
\begin{equation}
\label{Hel(,)}
H(\mu,\nu)=\int_{X}\sqrt{\frac{d\mu}{d\rho}\frac{d\nu}{d\rho}}d\rho=
\left(\frac{{\rm det}\,B_n{\rm det}\,C_n}{{\rm
det}^2\,\frac{B_n+C_n}{2}}\right)^{1/4}.
\end{equation}
The latter formula is based on the following formula for a positive definite operator $C$
in the space ${\mathbb R}^m$:
\begin{equation}
\frac{1}{\sqrt{\pi^m}}\int_{{\mathbb
R}^m}\exp(-(Cx,x))dx=\frac{1}{\sqrt{{\rm det}\,C}}.
\end{equation}
Let, as before, $t=(t_{rs})_{r,s=1}^m\in{\rm GL}(m,{\mathbb R}),\,\,
B_n={\rm
diag}(b_{1n},b_{2n},...,b_{mn}),\,\,X_n(t)\!=\!B_n^{1/2}tB_n^{-1/2}$
$\in {\rm Mat}(m,{\mathbb R})$.
Let $M^{i_1i_2...i_r}_{j_1j_2...j_r}(t)$ be the minors of the
matrix $t$ with $i_1,i_2,...,i_r$ rows and $j_1,j_2,...,j_r$
columns.

Let us denote by $\mu^{(B,a)}=\mu_{(C,a)}$  the Gaussian measure with
the covariance operator $C=(2B)^{-1}$ on the space ${\mathbb R}^m$
defined by the formula: $ \frac{d\mu^B(x)}{dx}= $
\begin{equation}
\label{dmu}
\sqrt{\frac{{\rm det}\,B}{\pi^m}}
\exp\Big(\!\!-\big(Bx,x\big)\Big)\!=\! \frac{1}{\sqrt{(2\pi)^m{\rm
det}\,C}}\exp\Big(\!\!-\frac{1}{2}\big(C^{-1}x,x)\big)\Big)\!=\!\frac{d\mu_{C}(x)}{dx}.
\end{equation}
Recall that by definition $\mu^f(\Delta)=\mu(f^{-1}(\Delta))$.
Since $L_tx=tx$, we get $\mu^{L^{-1}_t}(x)=\mu(tx)$ therefore,
\begin{equation}
\label{dmut}
\left(\mu_m^{(B_n,a)}\right)^{L^{-1}_t}(x)=\mu_m^{(B_n(t),t^{-1}a)}\quad\text{where}\quad
B_n(t)=t^*B_nt.
\end{equation}
Indeed,
$$
d\left(\mu_m^{(B_n,a)}\right)^{L^{-1}_t}(x)= \sqrt{\frac{{\rm
det}\,B_n}{\pi^m} }\exp\Big(-\big(B_nt(x-t^{-1}a),t(x-t^{-1}a)\big)\Big)dtx=
$$
$$
 \sqrt{\frac{\vert{\rm det}\,\,t\vert^2{\rm
det}\,B_n}{\pi^m} }\exp\Big(-\big(t^*B_nt(x-t^{-1}a),(x-t^{-1}a)\big)\Big)dx=
d\mu_m^{(B_n(t),t^{-1}a)}(x),
$$
where $B_n(t)=t^*B_nt,\,\,B_n={\rm
diag}(b_{1n},b_{2n},...,b_{mn}),\,\,{\rm det}\,B_n(t)=\vert{\rm
det}\,\,t\vert^2{\rm det}\,B_n.$

Using (\ref{Hel(,)}),  (\ref{dmu}) and (\ref{dmut}) we obtain
$$
H\left(\left(\mu_m^{(B_n,0)}\right)^{L_t},\mu_m^{(B_n,0)}\right)=H\left(\mu_m^{(B_n(t),0)},\mu_m^{(B_n,0)}\right)=
\left(\frac{{\rm det}\,B_n(t){\rm
det}\,B_n}{\pi^m\pi^m}\right)^{1/4} \times
$$
$$
\int_{{\mathbb
R}^m}\!\!\!\exp\Big(\!-\!\Big(\frac{B_n(t)\!+\! B_n}{2}x,x\Big)\Big)dx\!
=\!\left(\frac{{\rm det}\,B_n(t){\rm det}\,B_n}{{\rm
det}^2\,\frac{B_n(t)+B_n}{2}}\right)^{1/4}\!\!\!\!\!=\!\!\left(\frac{{\rm
det}\,C_n(t)}{\vert{\rm det}\,t\vert{\rm det}\,B_n}\right)^{-1/2}\!\!\!,
$$
where $C_n(t)=\frac{B_n(t)+B_n}{2}=\frac{t^*B_nt+B_n}{2}.$
Now we  show that
\begin{equation}
\frac{{\rm det}\,C_n(t)}{\vert{\rm det}\,t\vert{\rm det}\,B_n}=
\frac{1}{2^m\vert{\rm det}\,\,t\vert}\,{\rm
det}\left(I+X_n^*(t)X_n(t)\right),
\end{equation}
where $X_n(t)=B_n^{1/2}tB_n^{-1/2}$. The latter equation is
equivalent to
$$
\frac{{\rm det}\,(t^*B_nt+B_n)}{{\rm det}\,B_n}= {\rm
det}\left(I+X_n^*(t)X_n(t)\right).
$$
To complete the proof of (\ref{Hel1-sl}) it is sufficient to see
that
$$
I+X_n^*(t)X_n(t)=I+B_n^{-1/2}t^*B_n^{1/2}B_n^{1/2}tB_n^{-1/2}=
B_n^{-1/2}(B_n+t^*B_nt)B_n^{-1/2}.
$$
The proof of  relation (\ref{Hel2-sl}) is based on the following theorem
that one can find, e.g., in \cite[ Ch. III, \S 16, Theorem 2]{Sko74}. 
\begin{thm}
\label{t.(B,a)sim(B,0)}
 Two Gaussian measures $\mu_{B,a}$ and $\mu_{B,b}$ are
equivalent if and only of $B^{-1/2}(a-b)\in H$.
\end{thm}
{
\small Indeed, we have
$$
\Vert C^{-1/2}(ta-a)\Vert^2_H\!=\!\!\sum_{n\in {\mathbb Z}}
\Vert C_n^{-1/2}(t-I)a_n\Vert^2_{H_n}\!=\!2\!\sum_{n\in
{\mathbb Z}}\sum_{r=1}^m
\frac{b_{kn}}{d_{kn}}\Big(\sum_{s=1}^m(t_{rs}\!-\delta_{rs})a_{sn}\Big)^2\!d_{kn}.
$$
To explain the latter equality let us describe $H$ and $C$. To
find an operator $C$ we present the measure $\mu^m_{(b,a)}$ in the
canonical form $\mu_{C,a}$ defined by its Fourier transform:
\begin{equation}
\label{A.F(H,mu_B),}
 \int_{H}\exp i(y,x)d\mu_{C,a}(x)=\exp\left(i(a,y)-\frac{1}{2}(Cy,y)\right),\,\,
y\in H,
\end{equation}
}
where $C$ is a positive {\it nuclear operator} (called the {\it
covariance operator}) on the Hilbert space $H$, and $a\in H$ is
the {\it mathematical expectation} or {\it mean}.

Recall the {\it Kolmogorov zero-one law}. Let us consider in the space ${\mathbb R}^\infty={\mathbb
R}\times {\mathbb R}\times\cdots$ the infinite tensor product  $\mu_b=\otimes_{n\in
{\mathbb N}}\mu_{b_k}$  of one-dimensional Gaussian measures $\mu_{b_k}$ on $\mathbb R$ defined as follows:
\begin{equation}
\label{A.mu-(b,a)}
d\mu_{b}(x)=\sqrt{b/\pi}\exp(-bx^2)dx.
\end{equation}
Consider a Hilbert space $l_2(a)$ defined by
\begin{equation*}
l_2(a)=\big\{x\in {\mathbb R}^\infty\,:\, \Vert x
\Vert^2_{l_2(a)}=\sum_{k\in {\mathbb N}}x_k^2a_k<\infty\big\},
\end{equation*}
where $a=(a_k)_{k\in {\mathbb N}}$ is an infinite sequence of
positive numbers.
\begin{thm}[Kolmogorov's zero-one law, \cite{ShFDT67}]
\label{A.Kol-0-1}
We have
$$
\mu_b(l_2(a))= \left\{\begin{array}{ccc} 0,&\text{if}&\sum_{k\in
{\mathbb N}}\frac{a_k}{b_k}=\infty,\\
1,&\text{if}&\sum_{k\in {\mathbb N}}\frac{a_k}{b_k}<\infty.
\end{array}\right.
$$
\end{thm}
Define the Hilbert
space $H\subset X_m$ as follows:
$$
H=l_2({\mathbb R}^m,d)=\big\{x=(x_{kn})_{k,n}\in X_m\mid\Vert x
\Vert^2_H:=\sum_{1\leq k\leq m,n\in {\mathbb
N}}x_{kn}^2d_{kn}<\infty\big\},
$$
where a sequence $d=(d_{kn})_{1\leq k\leq
m,n\in{\mathbb Z}}$ of positive numbers is chosen such that $ \sum_{1\leq
k\leq m,n\in {\mathbb N}}\frac{d_{kn}}{b_{kn}}<\infty$. Then by
the Kolmogorov zero-one law,
$\mu^m_{(b,a)}(H)=1$.
We show that $C={\rm diag}(c_{kn})$, where
$c_{kn}=\frac{d_{kn}}{2b_{kn}}$. Indeed, we get
$$
\sum_{1\leq k\leq m,n\in {\mathbb N}}b_{kn}x_{kn}^2=
\frac{1}{2}\sum_{1\leq k\leq m,n\in {\mathbb N}}\frac{2b_{kn}}{d_{kn}}x_{kn}^2d_{kn}=\frac{1}{2}(C^{-1}x,x)_H.
$$

{\it Proof} of  Lemma \ref{ldet*}. Let as recall the definition of
the Gram  determinant and the Gram  matrix   (see
\cite{Gan58}, Chap IX, \S 5). For vectors $x_1,x_2,..., x_m$ in some
Hilbert space $H$ the Gram
 matrix $\gamma(x_1,x_2,..., x_m)$ is defined  by the formula
$$
\gamma(x_1,x_2,..., x_m)=((x_k,x_n)_{k,n=1}^m).
$$
\index{Gram determinant}\index{Gram matrix}
The determinant of this matrix is called the Gram  determinant for
the vectors $x_1,x_2,..., x_m$ and is denoted by
$\Gamma(x_1,x_2,..., x_m)$. Thus,
$$
\Gamma(x_1,x_2,..., x_m):={\rm det}\,\gamma(x_1,x_2,..., x_m).
$$
Let 
$$
X= \left(
\begin{array}{cccc}
x_{11}&x_{12}&...&x_{1m}\\
x_{21}&x_{22}&...&x_{2m}\\
...   &...   &...&...\\
x_{m1}&x_{m2}&...&x_{mm}
\end{array}
\right).
$$
Set $x_k=(x_{1k},x_{2k},...,x_{mk})\in{\mathbb R}^m,\,1\leq k\leq
m,$ then, obviously, we get
$$
X^*X= \left(
\begin{array}{cccc}
(x_1,x_1)&(x_1,x_2)&...&(x_1,x_m)\\
(x_2,x_1)&(x_2,x_2)&...&(x_2,x_m)\\
...   &...   &...&...\\
(x_m,x_1)&(x_m,x_2)&...&(x_m,x_m)
\end{array}
\right)=\gamma(x_1,x_2,..., x_m).
$$
We would like to find an exact expression for ${\rm
det}\left(I+\gamma(x_1,x_2,..., x_m)\right).$ It is convenient to
consider the following function:
$$
F_{m,X}^\lambda=
F^{\,\,\,\lambda_1,\lambda_2,...,\lambda_m}_{m;x_1,x_2,...,x_m}=
{\rm det}\Big(\sum_{k=1}^m\lambda_kE_{kk}+\gamma(x_1,x_2,...,
x_m)\Big),\quad \lambda\in {\mathbb C}^m.
$$
It is easy to see that for $m=2$ we have
$$
F_{2;x_1,x_2}^{\,\,\,\lambda_1,\lambda_2}=
{\rm det}\left(
\begin{array}{cc}
\lambda_1+(x_1,x_1)&(x_1,x_2)\\
(x_2,x_1)&\lambda_2+(x_2,x_2)
\end{array}
\right)=
$$
$$
\lambda_1\lambda_2+\lambda_1\Gamma(x_2)+\lambda_2\Gamma(x_1)+\Gamma(x_1,x_2)=
$$
\begin{equation}
\label{5.F_2}
\lambda_1\lambda_2\left(1+\lambda_1^{-1}\Gamma(x_1)+\lambda_2^{-1}\Gamma(x_2)+
(\lambda_1\lambda_2)^{-1}\Gamma(x_1,x_2)\right).
\end{equation}
The general formula is
\begin{equation}
\label{F_m(la,x)}
F^{\,\,\,\lambda_1,\lambda_2,...,\lambda_m}_{m;x_1,x_2,...,x_m}={\rm
det}\Big(\sum_{k=1}^m\lambda_kE_{kk}+\gamma(x_1,x_2,...,
x_m)\Big)=
\end{equation}
$$
\prod_{k=1}^m\lambda_k\Big(1+\sum_{r=1}^m\sum_{1\leq
i_1<i_2<...<i_r\leq
m}\Big(\lambda_{i_1}\lambda_{i_2}...\lambda_{i_r}\Big)^{-1}
\Gamma(x_{i_1},x_{i_2},...,x_{i_r})\Big)=
$$
$$
\prod_{k=1}^m\lambda_k\Big(1+\sum_{r=1}^{m} \sum_{1\leq
i_1<i_2<...<i_r\leq m;1\leq j_1<j_2<...<j_r\leq m}
\Big(\lambda_{i_1}\lambda_{i_2}...\lambda_{i_r}\Big)^{-1}
\Big(M^{i_1i_2...i_r}_{j_1j_2...j_r}(X)\Big)^2\Big).
$$
We have used the following formula (see \cite{Gan58}, Chap IX, \S
5 formula (25)):
\begin{equation}
\label{Gramm(x,y)=M^2(X)}
\Gamma(x_{i_1},x_{i_2},...,x_{i_r})= \sum_{1\leq
j_1<j_2<...<j_r\leq m}
\left(M^{i_1i_2...i_r}_{j_1j_2...j_r}(X)\right)^2.
\end{equation}
Finally, using (\ref{F_m(la,x)}) for
$(\lambda_1,\lambda_2,...,\lambda_m)=(1,1,...,1)$ we  get
(\ref{det*}).

We study the case $m=2$ more carefully.
\begin{lem}
\label{l7.5}
 For $t\in GL(2,{\mathbb R})$ we have, if ${\rm det}\,t>0$,
$$
(\mu_{(b,0)}^2)^{L_t}\perp\mu_{(b,0)}^2\quad\Leftrightarrow\quad
$$
\begin{equation}
\sum_{n\in{\mathbb Z}}\Big[ (1-\mid{\rm
det}\,t\mid)^2+(t_{11}-t_{22})^2+\Big(t_{12}\sqrt{\frac{b_{1n}}{b_{2n}}}+
t_{21}\sqrt{\frac{b_{2n}}{b_{1n}}}\Big)^2\Big]=\infty.
\end{equation}
If ${\rm det}\,t<0$ we have
$$
(\mu_{(b,0)}^2)^{L_t}\perp\mu_{(b,0)}^2\quad\Leftrightarrow\quad
$$
\begin{equation}
\sum_{n\in{\mathbb Z}} \Big[ (1-\mid{\rm det}\,t\mid)^2 +
(t_{11}+t_{22})^2+\Big(t_{12}\sqrt{\frac{b_{1n}}{b_{2n}}}-
t_{21}\sqrt{\frac{b_{2n}}{b_{1n}}}\Big)^2\Big]=\infty.
\end{equation}
\end{lem}
\begin{pf}   Using (\ref{Hel1-sl}) set
 $$
H_{m,n}(t)=H\left(\left(\mu_m^{(B_n,0)}\right)^{L^{-1}_t},\mu_m^{(B_n,0)}\right)=
\left(\frac{1}{2^m\vert{\rm det}\,\,t\vert} {\rm
det}\left(I+X_n^*(t)X_n(t)\right)\right)^{-1/2}.
$$
For $m=2$  using (\ref{5.X_n(t)}) we get $X(t)=B^{1/2}tB^{-1/2}$
hence,
$$
X(t)=
 \left(\begin{array}{cc}
 b_{1n}&0\\
 0&b_{2n}
\end{array}\right)^{1/2}\!\!
 \left(\begin{array}{cc}
 t_{11}&t_{12}\\
 t_{21}&t_{22}
\end{array}\right)
 \left(\begin{array}{cc}
 b_{1n}&0\\
 0&b_{2n}
\end{array}\right)^{-1/2}=\!\!
 \left(\begin{array}{cc}
 t_{11}&\sqrt{ \frac{b_{1n}}{ b_{2n}}}t_{12}\\
 \sqrt{\frac{ b_{2n}}{ b_{1n}}}t_{21}&t_{22}
\end{array}\right).
$$
Therefore, using  (\ref{5.F_2}) we get
$$
H_{2,n}^{-2}(t)=\frac{1}{2^2\mid{\rm det}\,t\mid}\left(1+\mid{\rm
det}\,t\mid^2+t_{11}^2+t_{22}^2+\frac{b_{1n}}{b_{2n}}t_{12}^2+
\frac{b_{2n}}{b_{1n}}t_{21}^2\right).
$$
Using Lemma \ref{l7.2} it is sufficient to calculate
$H_{2,n}^{-2}(t)-1.$ Indeed, for ${\rm det}\,t>0$ we have
$$
H_{2,n}^{-2}(t)-1=\frac{1}{2^2\mid{\rm det}\,t\mid}\times
$$
$$
\left( 1-2{\rm det}\,t+\mid{\rm
det}\,t\mid^2+t_{11}^2+t_{22}^2+\frac{b_{1n}}{b_{2n}}t_{12}^2+
\frac{b_{2n}}{b_{1n}}t_{21}^2-2(t_{11}t_{22}-t_{12}t_{21})\right)=
$$
$$
\frac{1}{2^2\mid{\rm det}\,t\mid} \Big[(1-\mid{\rm
det}\,t\mid)^2+
(t_{11}-t_{22})^2+\Big(t_{12}\sqrt{\frac{b_{1n}}{b_{2n}}}+
t_{21}\sqrt{\frac{b_{2n}}{b_{1n}}}\Big)^2\Big].
$$
For ${\rm det}\,t<0$ we get
$$
H_{2,n}^{-2}(t)-1=\frac{1}{2^2\mid{\rm det}\,t\mid}\times
$$
$$
\Big( 1+2{\rm det}\,t+\mid{\rm
det}\,t\mid^2+t_{11}^2+t_{22}^2+\frac{b_{1n}}{b_{2n}}t_{12}^2+
\frac{b_{2n}}{b_{1n}}t_{21}^2+2(t_{11}t_{22}-t_{12}t_{21})\Big)=
$$
$$
\frac{1}{2^2\mid{\rm det}\,t\mid} \Big[(1-\mid{\rm
det}\,t\mid)^2+
(t_{11}+t_{22})^2+\Big(t_{12}\sqrt{\frac{b_{1n}}{b_{2n}}}-
t_{21}\sqrt{\frac{b_{2n}}{b_{1n}}}\Big)^2\Big].
$$
\qed\end{pf}
Using Lemma \ref{l7.3}, Lemma \ref{l7.5} and (\ref{Hel2-sl}) we
get
\begin{lem}
\label{l7.6}
 For $t\in GL(2,{\mathbb R})$ we have
\begin{equation*}
\label{....}
 (\mu_{(b,a)}^2)^{L_t}\perp\mu_{(b,a)}^2\quad{\it
if}\quad\mid{\rm det}\,t\mid\not=1.
\end{equation*}
If ${\rm det}\,t=1$, we have
$$
(\mu_{(b,a)}^2)^{L_t}\perp\mu_{(b,a)}^2\quad\Leftrightarrow\quad\Sigma^+(t)=\Sigma_1^+(t)+
\Sigma_2(t)=\infty.
$$
If ${\rm det}\,t=-1$, we have
$$
(\mu_{(b,a)}^2)^{L_t}\perp\mu_{(b,a)}^2\quad\Leftrightarrow\quad\Sigma^{-}(t)=\Sigma_1^{-}(t)+
\Sigma_2(t)=\infty,
$$
where
$$
\Sigma_1^+(t)=\sum_{n\in{\mathbb
Z}}\Big[(t_{11}-t_{22})^2+\Big(t_{12}\sqrt{\frac{b_{1n}}{b_{2n}}}+
t_{21}\sqrt{\frac{b_{2n}}{b_{1n}}}\Big)^2\Big],
$$
$$
\Sigma_1^{-}(t)=\sum_{n\in{\mathbb
Z}}\Big[(t_{11}+t_{22})^2+\Big(t_{12}\sqrt{\frac{b_{1n}}{b_{2n}}}-
t_{21}\sqrt{\frac{b_{2n}}{b_{1n}}}\Big)^2\Big],
$$
\begin{equation}
\label{Sigma_2}
\Sigma_2(t^{-1})=\sum_{n\in{\mathbb
Z}}\Big[b_{1n}\big[(t_{11}-1)a_{1n}+t_{12}a_{2n}\big]^2+b_{2n}\big[t_{21}a_{1n}+(t_{22}-1)a_{2n}\big]^2\Big].
\end{equation}
\end{lem}
%
\begin{pf} {\it of  Lemma \ref{perp2}}.
We show that it is sufficient to consider only five
particular cases:
$$
\exp(tE_{12})=I+tE_{12}= \left(\begin{array}{cc}
 1&t\\
 0&1
\end{array}\right),\quad
\exp(tE_{21})=I+tE_{21}=\left(\begin{array}{cc}
 1&0\\
 t&1
\end{array}\right),
$$
$$
\exp(tE_{12})P_1= \left(\begin{array}{cc}
 -1&t\\
  0&1
\end{array}\right),\quad
\exp(tE_{21})P_2= \left(\begin{array}{cc}
 1&0\\
 t&-1
\end{array}\right),
$$
and
$$
\tau_{-}(\phi,s)= \left(
\begin{array}{cc}
\cos\phi&s^2\sin\phi\\
s^{-2}\sin\phi&-\cos\phi
\end{array}
\right),
$$
where
$$
P_1= \left(\begin{array}{cc}
 -1&0\\
  0&1
\end{array}\right),\quad
P_2=
 \left(\begin{array}{cc}
 1& 0\\
 0&-1
\end{array}\right).
$$
We note that $\tau_{-}(\phi,s)=$
$$
 \left(
\begin{array}{cc}
\cos\phi&s^2\sin\phi\\
s^{-2}\sin\phi&-\cos\phi
\end{array}
\right)= \left(\begin{array}{cc}
s&0\\
0&s^{-1}
\end{array}\right)
\left(\begin{array}{cc}
\cos\phi&-\sin\phi\\
\sin\phi&\cos\phi
\end{array}
\right) \left(
\begin{array}{cc}
s^{-1}&0\\
0     &s
\end{array}
\right)P_2.
$$
Using Lemma \ref{l7.5} we see that we have to consider only two
special cases:
$$
 t\in GL(2,{\mathbb R}),\quad {\rm det}\,t=1,\quad
t_{11}=t_{22},
$$ 
and
$$
 t\in GL(2,{\mathbb R}),\quad {\rm det}\,t=-1,\quad
t_{11}=-t_{22}.
$$
In the first case we have
$$
t=\left(
\begin{array}{cc}
\alpha     &t_{12}\\
t_{21}&\alpha
\end{array}
\right),\quad{\rm det}\,t=\alpha ^2-t_{12}t_{21}=1.
$$
In the second case we have
$$
t=\left(
\begin{array}{cc}
\alpha &t_{12}\\
t_{21}&-\alpha
\end{array}
\right),\quad{\rm det}\,t=-\alpha ^2-t_{12}t_{21}=-1.
$$
We can see that in the first (respectively second) case,  when
$t_{12}t_{21}\!>0$ (respectively $t_{12}t_{21}<0$), we have
$\Sigma_1^+(t)=\infty$ (respectively $\Sigma_1^{-}(t)=\infty$).

Indeed, if  ${\rm det}\,t=1$ and $t_{12}t_{21}\geq 1$,  then
$|t_{21}|\geq  | t_{12}|^{-1}$ and we have
$$
\left| t_{12}\sqrt{\frac{b_{1n}}{b_{2n}}}+
t_{21}\sqrt{\frac{b_{2n}}{b_{1n}}}\right|=
|t_{12}|\sqrt{\frac{b_{1n}}{b_{2n}}}+
|t_{21}|\sqrt{\frac{b_{2n}}{b_{1n}}}\geq
|t_{12}|\sqrt{\frac{b_{1n}}{b_{2n}}}+ |
t_{12}|^{-1}\sqrt{\frac{b_{2n}}{b_{1n}}}\geq 2.
$$
When  ${\rm det}\,t=1$ and $t_{12}t_{21}\in (0,1)$, then
$|t_{12}|^{-1}>|t_{21}|$ and we get
$$
\left| t_{12}\sqrt{\frac{b_{1n}}{b_{2n}}}+
t_{21}\sqrt{\frac{b_{2n}}{b_{1n}}}\right|=t_{12}t_{21}\left(
|t_{21}|^{-1}\sqrt{\frac{b_{1n}}{b_{2n}}}+
|t_{12}|^{-1}\sqrt{\frac{b_{2n}}{b_{1n}}}\right)\geq 2|t_{12}t_{21}|.
$$
The same is true for the second case, i.e., when ${\rm det}\,t=-1$
and $t_{12}t_{21}<0$.

When
$$
{\rm det}\,t=\alpha^2-t_{12}t_{21}=1,\quad\text{and}\quad t_{12}t_{21}=0,
$$
we have four cases
\begin{equation}
\label{+1-4}
 \left(\begin{array}{cc}
 1& t\\
 0&1
\end{array}\right),\quad
\left(\begin{array}{cc}
 1&0\\
 t&1
\end{array}\right),\quad
\left(\begin{array}{cc}
 -1& t\\
 0 &-1
\end{array}\right),\quad
\left(\begin{array}{cc}
 -1& 0\\
 t &-1
\end{array}\right),\quad t\in {\mathbb R}.
\end{equation}
 When
$$
{\rm det}\,t=-\alpha^2-t_{12}t_{21}=-1,\quad\text{and}\quad t_{12}t_{21}=0,
$$
we also have four cases:
\begin{equation}
\label{-1-4}
\left(\begin{array}{cc}
 -1& t\\
 0 &1
\end{array}\right),\quad
\left(\begin{array}{cc}
 1&0\\
 t&-1
\end{array}\right),\quad
 \left(\begin{array}{cc}
 1& t\\
 0&-1
\end{array}\right),\quad
\left(\begin{array}{cc}
 -1& 0\\
 t &1
\end{array}\right),
\quad t\in {\mathbb R}.
\end{equation}
Thus, it remains to consider two cases:
$$
{\rm det}\,t=\alpha ^2-t_{12}t_{21}=1,\quad\text{and}\quad t_{12}t_{21}\in[-1,\,0),
$$
$$
{\rm det}\,t=-\alpha ^2-t_{12}t_{21}=-1,\quad\text{and}\quad t_{12}t_{21}\in (0,\,1].
$$
Finally, we can set in the first case $\alpha =\cos\phi$ since
$\alpha ^2=1+t_{12}t_{21}\in [0,1)$. Then
$-t_{12}t_{21}=\sin^2\phi$ so, $t_{12}=-s^2\sin\phi$ and
$t_{21}=s^{-2}\sin\phi$, with $s>0$.

In the second case we can set $\alpha =\cos\phi$ since
$\alpha^2=1-t_{12}t_{21}\in [0,1)$. Then $t_{12}t_{21}=\sin^2\phi$
so $t_{12}=s^2\sin\phi$ and $t_{21}=s^{-2}\sin\phi$, with $s>0$.
Finally, in the first (the second) case we have to consider
\begin{equation}
\label{tau+-()}
t\!=\!\tau_{+}(\phi,s)\!=\!\left(
\begin{array}{cc}
\cos\phi&-s^{2}\sin\phi\\
s^{-2}\sin\phi&\cos\phi
\end{array}
\right),\,\,\,
 t\!=\!\tau_{-}(\phi,s)=\!\left(
\begin{array}{cc}
\cos\phi&s^2\sin\phi\\
s^{-2}\sin\phi&-\cos\phi
\end{array}
\right).
\end{equation}
We show that only the first {\it two cases in (\ref{+1-4}) and
(\ref{-1-4}) and the  second case in (\ref{tau+-()}) are
independent}. Indeed, we have for $a=\left(\begin{array}{c}
 a_{1n}\\
 a_{2n}
\end{array}\right)$ (see Lemma~\ref{l7.6})
$$
t=\left(\begin{array}{cc}
 1& t\\
 0&1
\end{array}\right),\,\,
(t^{-1}-I)a=
 \left(\left(\begin{array}{cc}
 1& -t\\
 0&1
\end{array}\right)-I\right)
\left(\begin{array}{c}
 a_{1n}\\
 a_{2n}
\end{array}\right)=
\left(\begin{array}{c}
 -ta_{2n}\\
 0
\end{array}\right),
$$
$$
t=\left(\begin{array}{cc}
 1& 0\\
 t&1
\end{array}\right),\quad
(t^{-1}-I)a=\left(\left(\begin{array}{cc}
 1&0\\
 -t&1
\end{array}\right)-I\right)\left(\begin{array}{c}
 a_{1n}\\
 a_{2n}
\end{array}\right)=
\left(\begin{array}{c}
 0\\
 -ta_{1n}
\end{array}\right),
$$
$$
t=\left(\begin{array}{cc}
 -1&t\\
 0&1
\end{array}\right)=\left(\begin{array}{cc}
 1&t\\
 0&1
\end{array}\right)\left(\begin{array}{cc}
 -1&0\\
 0&1
\end{array}\right)
,\quad t^{-1}=\left(\begin{array}{cc}
 -1&0\\
 0&1
\end{array}\right)\left(\begin{array}{cc}
 1&-t\\
 0&1
\end{array}\right)=
$$
$$
\left(\begin{array}{cc}
 -1&t\\
 0&1
\end{array}\right),\quad
(t^{-1}-I)a=\left(\begin{array}{cc}
-2&t\\
 0&0\\
\end{array}\right)\left(\begin{array}{c}
 a_{1n}\\
 a_{2n}
\end{array}\right)=
\left(\begin{array}{c}
-2a_{1n}+ta_{2n}\\
0
\end{array}\right),
$$
$$
t=\left(\begin{array}{cc}
 1&0\\
 t&-1
\end{array}\right)=\left(\begin{array}{cc}
 1&0\\
 t&1
\end{array}\right)\left(\begin{array}{cc}
 1&0\\
 0&-1
\end{array}\right),\quad t^{-1}=
\left(\begin{array}{cc}
 1&0\\
 0&-1
\end{array}\right)\left(\begin{array}{cc}
 1&0\\
 -t&1
\end{array}\right)=
$$
$$
\left(\begin{array}{cc}
 1&0\\
 t&-1
\end{array}\right),\quad (t^{-1}-I)a=\left(\begin{array}{cc}
0&0\\
 t&-2\\
\end{array}\right)\left(\begin{array}{c}
 a_{1n}\\
 a_{2n}
\end{array}\right)=
\left(\begin{array}{c}
0\\
ta_{1n}-2a_{2n}
\end{array}\right).
$$
Therefore, we get
$$
\Sigma^+\left(\begin{array}{cc}
 1& t\\
 0&1
\end{array}\right)=
t^2\sum_{n\in{\mathbb
Z}}b_{1n}\left(\frac{1}{b_{2n}}+a_{2n}^2\right)\simeq
 S^L_{12}(\mu)= \sum_{n\in{\mathbb
Z}}\frac{b_{1n}}{2}\left(\frac{1}{2b_{2n}}+a_{2n}^2\right),\,\,t\not=0,
$$
$$
\Sigma^+\left(\begin{array}{cc}
 1& 0\\
 t&1
\end{array}\right)=
t^2\sum_{n\in{\mathbb
Z}}b_{2n}\left(\frac{1}{b_{1n}}+a_{1n}^2\right)\simeq
S^L_{21}(\mu)=\sum_{n\in{\mathbb
Z}}\frac{b_{2n}}{2}\left(\frac{1}{2b_{1n}}+a_{1n}^2\right),\,\,t\not=0,
$$
$$
\Sigma^{-}\left(\begin{array}{cc}
 -1& t\\
 0&1
\end{array}\right)
=t^2\sum_{n\in{\mathbb Z}}\frac{b_{1n}}{b_{2n}}+
\sum_{n\in{\mathbb Z}} b_{1n}(-2a_{1n}+ta_{2n})^2
=:S^{L,-}_{12}(\mu,t),
$$
$$
\Sigma^{-}\left(\begin{array}{cc}
 1& 0\\
 t&-1
\end{array}\right)
=t^2\sum_{n\in{\mathbb Z}}\frac{b_{2n}}{b_{1n}}+
\sum_{n\in{\mathbb Z}} b_{2n}(ta_{1n}-2a_{2n})^2
=:S^{L,-}_{21}(\mu,t).
$$
For the last two cases in (\ref{+1-4}) and (\ref{-1-4}) we get respectively
$$
t=\left(\begin{array}{cc}
 -1&t\\
 0&-1
\end{array}\right)=\left(\begin{array}{cc}
 -1&0\\
 0&-1
\end{array}\right)\left(\begin{array}{cc}
 1&-t\\
 0&1
\end{array}\right),\quad t^{-1}=\left(\begin{array}{cc}
 1&t\\
 0&1
\end{array}\right)\left(\begin{array}{cc}
 -1&0\\
 0&-1
\end{array}\right)=
$$
$$
\left(\begin{array}{cc}
 -1&-t\\
 0&-1
\end{array}\right),\quad
(t^{-1}-I)a=\left(\begin{array}{cc}
 -2&-t\\
 0&-2
\end{array}\right)\left(\begin{array}{c}
 a_{1n}\\
 a_{2n}
\end{array}\right)=-\left(\begin{array}{c}
 2a_{1n}+t a_{2n}\\
2a_{2n}
\end{array}\right),
$$
$$
 t=\left(\begin{array}{cc}
 -1&0\\
 t&-1
\end{array}\right)=\left(\begin{array}{cc}
 -1&0\\
 0&-1
\end{array}\right)\left(\begin{array}{cc}
 1&0\\
 -t&1
\end{array}\right),\quad t^{-1}=\left(\begin{array}{cc}
 1&0\\
 t&1
\end{array}\right)\left(\begin{array}{cc}
 -1&0\\
 0&-1
\end{array}\right)
$$
$$
=\left(\begin{array}{cc}
 -1&0\\
 -t&-1
\end{array}\right),\quad (t^{-1}-I)a=\left(\begin{array}{cc}
-2&0\\
 -t&-2\\
\end{array}\right)\left(\begin{array}{c}
 a_{1n}\\
 a_{2n}
\end{array}\right)=
-\left(\begin{array}{c}
2a_{1n}\\
ta_{1n}+2a_{2n}
\end{array}\right).
$$
$$
t=\left(\begin{array}{cc}
 1&t\\
 0&-1
\end{array}\right)=\left(\begin{array}{cc}
 1&0\\
 0&-1
\end{array}\right)\left(\begin{array}{cc}
 1&t\\
 0&1
\end{array}\right),\quad t^{-1}=\left(\begin{array}{cc}
 1&-t\\
 0&1
\end{array}\right)\left(\begin{array}{cc}
 1&0\\
 0&-1
\end{array}\right)=
$$
$$
\left(\begin{array}{cc}
 1&t\\
 0&-1
\end{array}\right),\quad (t^{-1}-I)a=\left(\begin{array}{cc}
 0&t\\
 0&-2
\end{array}\right)\left(\begin{array}{c}
 a_{1n}\\
 a_{2n}
\end{array}\right)=
\left(\begin{array}{c}
ta_{2n}\\
-2a_{2n}
\end{array}\right),
$$
$$
t=\left(\begin{array}{cc}
 -1&0\\
 t&1
\end{array}\right)=\left(\begin{array}{cc}
 -1&0\\
 0&1
\end{array}\right)\left(\begin{array}{cc}
 1&0\\
 t&1
\end{array}\right),\quad t^{-1}=\left(\begin{array}{cc}
 1&0\\
 -t&1
\end{array}\right)\left(\begin{array}{cc}
 -1&0\\
 0&1
\end{array}\right)=
$$
$$
\left(\begin{array}{cc}
 -1&0\\
 t&1
\end{array}\right),\quad (t^{-1}-I)a=\left(\begin{array}{cc}
-2&0\\
 t&0\\
\end{array}\right)\left(\begin{array}{c}
 a_{1n}\\
 a_{2n}
\end{array}\right)=
\left(\begin{array}{c}
-2a_{1n}\\
ta_{1n}
\end{array}\right).
$$
Set
\begin{equation}
\label{S^L(11)-(22)}
 S^L_{11}(\mu):=S^{L,-}_{12}(\mu,0)=4\sum_{n\in{\mathbb
Z}}b_{1n}a_{1n}^2,\quad
S^L_{22}(\mu):=S^{L,-}_{21}(\mu,0)=4\sum_{n\in{\mathbb
Z}}b_{2n}a_{2n}^2.
\end{equation}
With this notation we see that the second two cases in
(\ref{+1-4}) and (\ref{-1-4}) are dependent:
$$
\Sigma^+\left(\begin{array}{cc}
 -1& t\\
 0&-1
\end{array}\right)
=t^2\sum_{n\in{\mathbb Z}}\frac{b_{1n}}{b_{2n}}+
\sum_{n\in{\mathbb
Z}}\left[b_{1n}(-2a_{1n}-ta_{2n})^2+b_{2n}(-2a_{2n})^2\right]
$$
$$
=S^{L,-}_{12}(\mu,-t)+S^L_{22}(\mu),\quad\text{note that}\quad \left(\begin{array}{cc}
 -1& t\\
 0 &-1
\end{array}\right)= \left(\begin{array}{cc}
 1& 0\\
 0 &-1
\end{array}\right)
 \left(\begin{array}{cc}
 -1& t\\
 0 &1
\end{array}\right).
$$
$$
\Sigma^+\left(\begin{array}{cc}
 -1& 0\\
 t&-1
\end{array}\right)
=t^2\sum_{n\in{\mathbb Z}}\frac{b_{2n}}{b_{1n}}+
\sum_{n\in{\mathbb
Z}}\left[b_{1n}(-2a_{1n})^2+b_{2n}(-ta_{1n}-2a_{2n})^2\right]
$$
$$
=S^{L,-}_{21}(\mu,-t)+S^L_{11}(\mu),\quad\text{note that}\quad
\left(\begin{array}{cc}
 -1& 0\\
 t &-1
\end{array}\right)=
\left(\begin{array}{cc}
 -1& 0\\
 0 &1
\end{array}\right)
\left(\begin{array}{cc}
 1& 0\\
 t &-1
\end{array}\right).
$$
$$
\Sigma^{-}\left(\begin{array}{cc}
 1& t\\
 0&-1
\end{array}\right)
=t^2\sum_{n\in{\mathbb Z}}\frac{b_{1n}}{b_{2n}}+
t^2\sum_{n\in{\mathbb Z}}b_{1n}a_{2n}^2+4\sum_{n\in{\mathbb Z}}b_{2n}a_{2n}^2\simeq t^2S^L_{12}(\mu)+S^L_{22}(\mu),
$$
$$
\Sigma^{-}\left(\begin{array}{cc}
 -1& 0\\
 t&1
\end{array}\right)
=t^2\sum_{n\in{\mathbb Z}}\frac{b_{2n}}{b_{1n}}+
4\sum_{n\in{\mathbb Z}} b_{1n}a_{1n}^2+t^2\sum_{n\in{\mathbb Z}} b_{2n}a_{1n}^2\simeq
t^2S^L_{21}(\mu)+S^L_{11}(\mu).
$$
To compare $(\mu_{(b,a)}^2)^{L_{\tau_{\pm}(\phi,s)}}$ and $\mu_{(b,a)}^2$
we calculate $\tau^{-1}_{+}(\phi,s)$ and $\tau^{-1}_{-}(\phi,s)$. Since
$$
\tau_{+}(\phi,s)=\left(
\begin{array}{cc}
\cos\phi&-s^{2}\sin\phi\\
s^{-2}\sin\phi&\cos\phi
\end{array}
\right),\quad \tau_{-}(\phi,s)= \left(
\begin{array}{cc}
\cos\phi&s^2\sin\phi\\
s^{-2}\sin\phi&-\cos\phi
\end{array}
\right),
$$
we get
$$
\tau^{-1}_{+}(\phi,s)=\left(
\begin{array}{cc}
\cos\phi&s^{2}\sin\phi\\
-s^{-2}\sin\phi&\cos\phi
\end{array}
\right),
\quad\tau^{-1}_{-}(\phi,s)=\left(
\begin{array}{cc}
\cos\phi&s^2\sin\phi\\
s^{-2}\sin\phi&-\cos\phi
\end{array}
\right)
$$
$=\tau_{-}(\phi,s)$.  Since $\tau^{-1}_{+}(\phi,s)-I=$
$$
 \left(
\begin{array}{cc}
\cos\phi-1&s^{2}\sin\phi\\
-s^{-2}\sin\phi&\cos\phi-1
\end{array}
\right)= \left(
\begin{array}{cc}
-2\sin^2\frac{\phi}{2}&s^{2}2\sin\frac{\phi}{2}\cos\frac{\phi}{2}\\
-s^{-2}2\sin\frac{\phi}{2}\cos\frac{\phi}{2}&-2\sin^2\frac{\phi}{2}
\end{array}
\right)=
$$
$$
\left(
\begin{array}{cc}
-2\sin\frac{\phi}{2}&0\\
0&-2\sin\frac{\phi}{2}
\end{array}
\right) \left(
\begin{array}{cc}
\sin\frac{\phi}{2}&-s^2\cos\frac{\phi}{2}\\
s^{-2}\cos\frac{\phi}{2}&\sin\frac{\phi}{2}
\end{array}
\right)
$$
and $ \tau^{-1}_{-}(\phi,s)-I=$
$$
 \left(
\begin{array}{cc}
\cos\phi-1&s^2\sin\phi\\
s^{-2}\sin\phi&-\cos\phi-1
\end{array}
\right)= \left(
\begin{array}{cc}
-2\sin^2\frac{\phi}{2}&s^{2}2\sin\frac{\phi}{2}\cos\frac{\phi}{2}\\
s^{-2}2\sin\frac{\phi}{2}\cos\frac{\phi}{2}&-2\cos^2\frac{\phi}{2}
\end{array}
\right)=
$$
$$
\left(
\begin{array}{cc}
-2\sin\frac{\phi}{2}&0\\
0&-2\cos\frac{\phi}{2}
\end{array}
\right) \left(
\begin{array}{cc}
\sin\frac{\phi}{2}&-s^2\cos\frac{\phi}{2}\\
-s^{-2}\sin\frac{\phi}{2}&\cos\frac{\phi}{2}
\end{array}
\right),
$$
we have (see (\ref{Sigma_2}))
$$ \Sigma_2(\tau_{-}(\phi,s))=
4\sin^2\frac{\phi}{2}\sum_{n\in {\mathbb Z}}
b_{1n}\Big(\sin\frac{\phi}{2}a_{1n}\!-\!s^2\cos\frac{\phi}{2}a_{2n}\Big)^2\!\!+
$$
$$
4\cos^2\frac{\phi}{2}\sum_{n\in {\mathbb Z}}b_{2n}\!
\Big(\!-s^{-2}\sin\frac{\phi}{2}a_{1n}+\cos\frac{\phi}{2}a_{2n}\Big)^2
$$
$$
\sim\sum_{n\in {\mathbb Z}}\big(4\sin^2\frac{\phi}{2}b_{1n}+4\cos^2\frac{\phi}{2}s^{-2}b_{2n}\big)
\Big(\sin\frac{\phi}{2}a_{1n}\!-\!s^2\cos\frac{\phi}{2}a_{2n}\Big)^2,
$$
\begin{equation}
\label{s(tau-)}
\Sigma_2(\tau_{-}(\phi,s))\!=\!\sum_{n\in {\mathbb Z}}\big(4\sin^2\frac{\phi}{2}b_{1n}\!+\!4\cos^2\frac{\phi}{2}s^{-2}b_{2n}\big)
\Big(\sin\frac{\phi}{2}a_{1n}\!-\!s^2\cos\frac{\phi}{2}a_{2n}\Big)^2.
\end{equation}
Finally, for $t=\tau_{-}(\phi,s)$ we get
\begin{equation}
\label{tau(phi,s)-}
\mu^{L_{\tau_-(\phi,s)}}\perp\mu \Leftrightarrow \sin^2\phi\Sigma_1(s)+\Sigma_2(\tau_{-}(\phi,s))=\infty,
\end{equation}
where
$$
\Sigma_1(s)=\sum_{n\in {\mathbb Z}}
\Big(s^2\sqrt{\frac{b_{1n}}{b_{2n}}}-s^{-2}\sqrt{\frac{b_{2n}}{b_{1n}}}\Big)^2.
$$
We have for $t=\tau_{+}(\phi,s)$ (see (\ref{Sigma_2}))
\begin{equation}
\label{tau(phi,s)+}
\mu^{L_{\tau_+(\phi,s)}}\perp\mu \Leftrightarrow \sin^2\phi\Sigma_1(s)
+\Sigma_2(\tau_{+}(\phi,s))=\infty,
\end{equation}
where $\Sigma_2(\tau_{+}(\phi,s))=$
\begin{equation}
\label{s(tau+)}
4\sin^2\frac{\phi}{2}\sum_{n\in {\mathbb Z}} \Big[
b_{1n}\Big(\sin\frac{\phi}{2}a_{1n}-s^2\cos\frac{\phi}{2}a_{2n}\Big)^2+
b_{2n}\Big(s^{-2}\cos\frac{\phi}{2}a_{1n}+\sin\frac{\phi}{2}a_{2n}\Big)^2
\Big].
\end{equation}
We show that the condition $\mu^{L_{\tau_+(\phi,s)}}\perp\mu$
depends on the previous conditions of the orthogonality. Indeed,
for $t=\tau_-(\phi,s)$ we have
$$
\mu^{L_{\tau_-(\phi,s)}}\perp\mu \Leftrightarrow (a)\,\, \Sigma_1(s)=\infty\quad\text{or}\quad(b)\,\,
\Sigma_1(s)<\infty,\,\,\,\text{but}\,\,\,\Sigma_2(\tau_{-}(\phi,s))=\infty.
$$
For $t=\tau_+(\phi,s)$ we get respectively
$$
\mu^{L_{\tau_+(\phi,s)}}\perp\mu \Leftrightarrow (c)\,\, \Sigma_1(s)=\infty\quad\text{or}\quad(d)\,\,
\Sigma_1(s)<\infty,\,\,\,\text{but}\,\,\,\Sigma_2(\tau_{+}(\phi,s))=\infty.
$$
We see that $(c)\Leftrightarrow(a)$. To investigate the condition
$(d)$ we observe that if $\Sigma_1(s)<\infty$, then
$\lim_{n\to\infty}s^2\sqrt{\frac{b_{1n}}{b_{2n}}}=1$ therefore,
we have  $b_{2n}\sim s^4b_{1n}$ hence, the following equivalence
holds: $\Sigma_2(\tau_{+}(\phi,s))=$
$$
4\sin^2\frac{\phi}{2}\sum_{n\in {\mathbb Z}} \Big[
b_{1n}\Big(\sin\frac{\phi}{2}a_{1n}-s^2\cos\frac{\phi}{2}a_{2n}\Big)^2+
b_{2n}\Big(s^{-2}\cos\frac{\phi}{2}a_{1n}+\sin\frac{\phi}{2}a_{2n}\Big)^2
\Big]
$$
$$
\sim 4\sin^2\frac{\phi}{2}\sum_{n\in {\mathbb Z}} \Big[
b_{1n}\Big(\sin\frac{\phi}{2}a_{1n}-s^2\cos\frac{\phi}{2}a_{2n}\Big)^2+
b_{1n}\Big(\cos\frac{\phi}{2}a_{1n}+s^{2}\sin\frac{\phi}{2}a_{2n}\Big)^2
\Big]=
$$
$$
4\sin^2\frac{\phi}{2}\sum_{n\in {\mathbb Z}}
b_{1n}\left[a_{1n}^2+s^4 a_{2n}^2\right]\sim
4\sin^2\frac{\phi}{2}\sum_{n\in {\mathbb Z}}
\left(b_{1n}a_{1n}^2+b_{2n}a_{2n}^2\right)=
$$
$$
\sin^2\frac{\phi}{2}\left[S^L_{11}(\mu)+S^L_{22}(\mu)\right].
$$

We see that condition $(d)$
follows from  the conditions
$S^L_{11}(\mu)\!\!=\!\!S^{L,-}_{12}(\mu,0)\!\!=\!\!\infty$ and
$S^L_{22}(\mu)=S^{L,-}_{21}(\mu,0)=\infty.$ \ This  completes the
proof of  Lemma \ref{perp2}. \qed\end{pf}
\subsection{The explicit expression for $(D^{-1}(\lambda)\mu,\mu)$}
\label{ch.5.7}
The following lemma will be systematically used in what follows.
\begin{lem} For the matrix $D(\lambda_1,\lambda_2,...,\lambda_m)$ defined below
\label{l7.7}
\begin{equation}
\label{D(lambda)}
D(\lambda_1,\lambda_2,...,\lambda_m)= \left(
\begin{array}{cccc}
1+\lambda_1&1          &...&1\\
1          &1+\lambda_2&...&1\\
&&...&\\
1          &          1&...&1+\lambda_m
\end{array}
\right)
\end{equation}
and $\mu=(\mu_k)_{k=1}^m\in {\mathbb R}^m$ we have
\begin{equation}
\label{Dmumu} (D^{-1}(\lambda_1,\lambda_2,...,\lambda_m)\mu,\mu)\!=\!
\frac{\sum_{k=1}^m\frac{\mu_k^2}{\lambda_k}+\sum_{1\leq k<n\leq
m}\frac{(\mu_k-\mu_n)^2}{\lambda_k\lambda_n}}{1+\sum_{k=1}^m\frac{1}{\lambda_k}}.
\end{equation}
\end{lem}
\begin{pf}
Let us set $d_m(\lambda_1,\lambda_2,...,\lambda_m)={\rm det
}\left(D(\lambda_1,\lambda_2,...,\lambda_m)\right).$
It is easy to
see that
\begin{equation}
\label{det(D)}
d_m(\lambda_1,\lambda_2,...,\lambda_m)=\prod_{k=1}^m\lambda_k\left(1+
\sum_{k=1}^m\frac{1}{\lambda_k}\right).
\end{equation}
For arbitrary $m$ we have
$$
D^{-1}(\lambda_1,\lambda_2,...,\lambda_m)= \left(
\begin{array}{cccc}
1+\lambda_1&1          &...&1\\
1          &1+\lambda_2&...&1\\
&&...&\\
1          &          1&...&1+\lambda_m
\end{array}
\right)^{-1}=\left(D^{-1}_{kn}\right)_{k,n=1}^m,
$$
where
$$
D^{-1}_{nn}=\frac{d_{m-1}(\lambda_1,...,\hat{\lambda}_n,...\lambda_m)}
{d_m(\lambda_1,\lambda_2,...,\lambda_m)}= \left(1+
\sum_{k=1}^m\frac{1}{\lambda_k}\right)^{-1}\frac{1}{\lambda_n}
\left(1+ \sum_{k=1,k\not=n}^m\frac{1}{\lambda_k}\right),
$$
$$
D^{-1}_{kn}=\frac{-d_{m-1}(\lambda_1,...,\hat{\lambda}_n,...\lambda_m)\vert_{\lambda_k=0}}
{d_m(\lambda_1,\lambda_2,...,\lambda_m)}=-\frac{1}{\lambda_k\lambda_n}\left(1+
\sum_{k=1}^m\frac{1}{\lambda_k}\right)^{-1},\,\,k\not=n,
$$
since using (\ref{det(D)}) we have
$$
d_{m-1}(\lambda_1,...,\hat{\lambda}_n,...\lambda_m)\vert_{\lambda_k=0}\!=\!
\lim_{\lambda_k\rightarrow
0}\prod_{p=1,p\not=n}^m\lambda_p\left(1+
\sum_{p=1,p\not=n}^m\frac{1}{\lambda_p}\right)=\frac{1}{\lambda_k\lambda_n}
\prod_{p=1}^m\lambda_p.
$$
Finally, we have for  $\mu=(\mu_1,\mu_2,...,\mu_m)\in {\mathbb
R}^m$
$$
(D^{-1}(\lambda_1,\lambda_2,...,\lambda_m)\mu,\mu)=
\sum_{k,n=1}^mD^{-1}_{kn}\mu_k\mu_n=
$$
$$
\left(1+ \sum_{k=1}^m\frac{1}{\lambda_k}\right)^{-1}\left[
\sum_{n=1}^m \frac{\mu_n^2}{\lambda_n} \left(1+
\sum_{k=1,k\not=n}^m\frac{1}{\lambda_k}\right)-2 \sum_{1\leq
k<n\leq m}\frac{\mu_k\mu_n}{\lambda_k\lambda_n}\right]=
$$
$$
\left(1+ \sum_{k=1}^m\frac{1}{\lambda_k}\right)^{-1} \left[
\sum_{n=1}^m \frac{\mu_n^2}{\lambda_n}+
 \sum_{1\leq k<n\leq
m}\frac{(\mu_k-\mu_n)^2}{\lambda_k\lambda_n} \right].
$$
\qed\end{pf}
\begin{rem}
\label{acci.obs.}
{\rm  Some useful observations}.
If we set $f_{(m)}=(f_k)_{k=1}^m$ and $g_{(m)}$ $=(g_k)_{k=1}^m$ where
$f_k=\frac{\mu_k}{\sqrt{\lambda_k}}$ and $g_k=\frac{1}{\sqrt{\lambda_k}}$ we can recognize that
$$
\sum_{n=1}^m\frac{\mu_n^2}{\lambda_n}=\Vert f_{(m)}\Vert^2=\Gamma(f_{(m)})
$$
and
$$
\sum_{1\leq k<n\leq
m}\frac{(\mu_k-\mu_n)^2}{\lambda_k\lambda_n}=\sum_{1\leq k<n\leq m}\left|
\begin{array}{cc}
f_k&f_n\\
g_k&g_n\\
\end{array}
\right|^2=\Gamma(f_{(m)},g_{(m)})
$$
since
$$
\left|
\begin{array}{cc}
f_k&f_n\\
g_k&g_n\\
\end{array}
\right|^2
=\left|
\begin{array}{cc}
\frac{\mu_k}{\sqrt{\lambda_k}}&\frac{\mu_n}{\sqrt{\lambda_n}}\\
\frac{1}{\sqrt{\lambda_k}}&\frac{1}{\sqrt{\lambda_n}}\\
\end{array}
\right|^2=\frac{(\mu_k-\mu_n)^2}{\lambda_k\lambda_n}.
$$
Set $\Delta(f,g)=\frac{\Gamma(f)+\Gamma(f,g)}{\Gamma(g)+1}$ for two vectors $f$ and $g$. Finally, we get
\begin{equation}
\label{acci.obs.D(l,m)=G/G}
(D^{-1}(\lambda_1,\lambda_2,...,\lambda_m)\mu,\mu)=\Delta(f_{(m)},g_{(m)}))=
\frac{\Gamma(f_{(m)})+\Gamma(f_{(m)},g_{(m)})}{\Gamma(g_{(m)})+1}.
\end{equation}
where
$\Gamma(f_1,f_2,\dots,f_n)$ is the Gram determinant and
$\gamma(f_1,f_2,\dots,f_n)$ is the Gram matrix of $n$ vectors
$f_1,f_2,\dots,f_n$ in a Hilbert space (see \cite{Gan58}).
\end{rem}
\index{Gramm determinant}\index{Gramm matrix}
\subsection{The proof of  Lemmas ~\ref{x1x1} -- \ref{x2,dA}}
\label{sec5.5}
\begin{pf} The proof of  Lemma \ref{d1} is based on    Lemma \ref{l7.7}.
We find out when the inclusion
$$
D_{1n}{\bf 1}\in\langle A_{kn}{\bf
1}=(x_{1k}D_{1n}+x_{2k}D_{2n}){\bf 1}\mid k\in {\mathbb Z}\rangle
$$
holds. Fix $m\in {\mathbb N}$, since $Mx_{1k}=a_{1k}$, we put $\sum_{k=-m}^mt_ka_{1k}=(t,b)=1$,
where $t=(t_k)_{k=-m}^m$ and $b=(a_{1k})_{k=-m}^m$. We have
$$
\Vert \big[\sum_{k=-m}^mt_k(x_{1k}D_{1n}+x_{2k}D_{2n})-D_{1n}\big]{\bf
1}\Vert^2=
$$
$$
\Vert \sum_{k=-m}^mt_k[(x_{1k}-a_{1k})D_{1n}+x_{2k}D_{2n}]{\bf
1}\Vert^2=\sum_{-m\leq k,r\leq m}(f_k,f_r)t_kt_r=:(A_{2m+1}t,t),
$$
where $A_{2m+1}=((f_k,f_r))_{k,r=-m}^m,$ and
$f_k=[(x_{1k}-a_{1k})D_{1n}+x_{2k}D_{2n}]{\bf 1}.$  We have
$$
(f_k,f_k)=\Vert\left[(x_{1k}-a_{1k})D_{1n}+x_{2k}D_{2n}\right]{\bf
1}\Vert^2=\frac{1}{2b_{1k}}\frac{b_{1n}}{2}+\Big(\frac{1}{2b_{2k}}+a_{2k}^2\Big)
\frac{b_{2n}}{2} \sim
$$
$$
\frac{1}{2b_{1k}}+\frac{1}{2b_{2k}}+a_{2k}^2,
$$
$$
(f_k,f_r)=\left(
\left[(x_{1k}-a_{1k})D_{1n}+x_{2k}D_{2n}\right]{\bf 1},
\left[(x_{1r}-a_{1r})D_{1n}+x_{2r}D_{2n}\right]{\bf 1}\right) =
$$
$$
(x_{2k},x_{2r})(D_{2n}{\bf 1},D_{2n}{\bf 1})=
a_{2k}a_{2r}\frac{b_{2n}}{2}\simeq a_{2k}a_{2r}.
$$
Finally, we have
\begin{equation}
\label{(A_{r,s}),D1,m=2}
(f_k,f_k)\sim\frac{1}{2b_{1k}}+\frac{1}{2b_{2k}}+a_{2k}^2,\quad
(f_k,f_r)\sim a_{2k}a_{2r},\quad k\not=r.
\end{equation}
For $A_{(m)}=((f_k,f_r))_{k,r=1}^m,$ and
$b=(a_{11},a_{12},...,a_{1m})\in {\mathbb R}^m$ we have
$$
A_{(m)}= \gamma(f_1,f_2,...,f_m)=\left(
\begin{array}{cccc}
(f_1,f_1)&(f_1,f_2)&...&(f_1,f_m)\\
(f_2,f_1)&(f_2,f_2)&...&(f_2,f_m)\\
             &             &...&             \\
(f_m,f_1)&(f_m,f_2)&...&(f_m,f_m)
\end{array}
\right)=
$$
$$
\left(
\begin{array}{cccc}
\frac{1}{2b_{11}}+\frac{1}{2b_{21}}+a_{21}^2& a_{21}a_{22}&...&a_{21}a_{2m}\\
 a_{22}a_{21}&\frac{1}{2b_{12}}+\frac{1}{2b_{22}}+a_{22}^2&...&a_{22}a_{2m}\\
                                            &             &...&             \\
 a_{2m}a_{21} & a_{2m}a_{22}&...& \frac{1}{2b_{1m}}+\frac{1}{2b_{2m}}+a_{2m}^2
\end{array}
\right)=
$$
$$
\left(\!\!
\begin{array}{cccc}
a_{21}&0&...&0\\
0&a_{22}&...&0\\
 &      &...& \\
0&0&...&a_{2m}
\end{array}
\!\!\right)
\!\!
 \left(\!\!
\begin{array}{cccc}
1+\lambda_1&1          &...&1\\
1          &1+\lambda_2&...&1\\
           &           &...& \\
1          &1&...&1+\lambda_m
\end{array}
\!\!\right)
\!\!
 \left(\!\!
\begin{array}{cccc}
a_{21}&0&...&0\\
0&a_{22}&...&0\\
 &      &...& \\
0&0&...&a_{2m}
\end{array}
\!\!\right),
$$
where
$\lambda_k=\frac{\frac{1}{2b_{1k}}+\frac{1}{2b_{2k}}}{a_{2k}^2},\,\,1\leq
k\leq m.$ Using  (\ref{D(lambda)}) we conclude that
$$
A_{(m)}={\rm
diag}(a_{21},a_{22},\dots,a_{2m})D(\lambda_1,\lambda_2,\dots,\lambda_m)
{\rm diag}(a_{21},a_{22},...,a_{2m}).
$$
Recall that $\mu={\rm
diag}(a_{21},a_{22},,...,a_{2m})^{-1}b=(\frac{a_{11}}{a_{21}},\frac{a_{12}}{a_{22}},...,
\frac{a_{1m}}{a_{2m}})$, where $b=$
\\$(a_{11},a_{12},...,a_{1m})\in{\mathbb R}^m$,  then
\begin{equation}
(A_{(m)}^{-1}b,b)\!=\!(D^{-1}(\lambda_1,\lambda_2,...,\lambda_m)\mu,\mu),\,\,
\lambda_k=\Big(\frac{1}{2b_{1k}}+\frac{1}{2b_{2k}}\Big)a_{2k}^{-2},\,\,\mu_k=a_{1k}a_{2k}^{-1}.
\end{equation}
Using Lemma \ref{l7.7} for  the  operator $A_{2m+1}$, and the vector $b\in {\mathbb R}^{2m+1}$ we obtain
$$
(A_{2m+1}^{-1}b,b)=\frac {\sum_{k=-m}^m\frac{a_{1k}^2}
{\frac{1}{2b_{1k}}+\frac{1}{2b_{2k}}}+\sum_{-m\leq k<n\leq m}
\frac{(a_{1k}a_{2n}-a_{1n}a_{2k})^2}
{\left(\frac{1}{2b_{1k}}+\frac{1}{2b_{2k}}\right)
\left(\frac{1}{2b_{1n}}+\frac{1}{2b_{2n}}\right)}} {
\sum_{k=-m}^m\frac{a_{2k}^2}
{\frac{1}{2b_{1k}}+\frac{1}{2b_{2k}}}+1 }
$$
$=\Delta(f_{(m)},g_{(m)})$, where
\begin{equation}
\label{f_m,g_m=1}
f_m\!=\!\Big(a_{1k}\Big(\frac{1}{2b_{1k}}+\frac{1}{2b_{2k}}\Big)^{-1/2}\Big)_{k=-m}^m,\,\,
g_m\!=\!\Big(a_{2k}\Big(\frac{1}{2b_{1k}}+\frac{1}{2b_{2k}}\Big)^{-1/2}\Big)_{k=-m}^m
\end{equation}
This proves  Lemma \ref{d1} \qed\end{pf}

The proof of  Lemma~\ref{d2} is exactly the same.

 The proof of  Lemma \ref{x1x1} is also based on  Lemma \ref{l7.7}. \par
\begin{pf} We study when
$
x_{1n}x_{1t}\in\langle A_{nk}A_{tk}{\bf 1}\mid k\in{\mathbb
Z}\rangle.
$
Since
$$
A_{nk}A_{tk}=(x_{1n}D_{1k}+x_{2n}D_{2k})(x_{1t}D_{1k}+x_{2t}D_{2k})=
$$
$$
x_{1n}x_{1t}D_{1k}^2+(x_{1n}x_{2t}+x_{2n}x_{1t})D_{1k}D_{2k}+x_{2n}x_{2t}D_{2k}^2
$$
and $MD_{1k}^2{\bf 1}=-\frac{b_{1k}}{2}$, set
$-\sum_{k=-m}^mt_k\frac{b_{1k}}{2}=(t,b')=1$, where $t=(t_k)_{k=-m}^m$ and $b'=-(\frac{b_{1k}}{2})_k\sim b=(b_{1k})_{k=-m}^m$. We
have
$$
\Vert \big[\sum_{k=-m}^mt_kA_{nk}A_{tk}-x_{1n}x_{1t}\big]{\bf
1}\Vert^2=
$$
$$
\Vert \sum_{k=-m}^mt_k\big[
x_{1n}x_{1t}\Big(D_{1k}^2+\frac{b_{1k}}{2}\big)+
(x_{1n}x_{2t}+x_{2n}x_{1t})D_{1k}D_{2k}+x_{2n}x_{2t}D_{2k}^2
\Big]{\bf 1}\Vert^2
$$
$$
=\sum_{-m\leq k,r\leq m}(f_k,f_r)t_kt_r=:(A_{2m+1}t,t),
$$
where $A_{2m+1}=((f_k,f_r))_{k,r=-m}^m$ and
$$
f_k= \Big[x_{1n}x_{1t}\Big(D_{1k}^2+\frac{b_{1k}}{2}\Big)+
(x_{1n}x_{2t}+x_{2n}x_{1t})D_{1k}D_{2k}+x_{2n}x_{2t}D_{2k}^2
\Big] {\bf 1} .
$$
 If we denote by $c_{kn}=\Vert x_{kn}\Vert^2=\frac{1}{2b_{kn}}+a_{kn}^2$,  we get
$$
(f_k,f_k)= \Vert\big[
x_{1n}x_{1t}\Big(D_{1k}^2+\frac{b_{1k}}{2}\Big)+
(x_{1n}x_{2t}+x_{2n}x_{1t})D_{1k}D_{2k}+x_{2n}x_{2t}D_{2k}^2
\big]{\bf 1}\Vert^2=
$$
$$
c_{1n}c_{1t}2\Big(\frac{b_{1k}}{2}\Big)^2 +
\big(c_{1n}c_{2t}+c_{1t}c_{2n}+2a_{1n}a_{2t}a_{1t}a_{2n}\big)\frac{b_{1k}}{2}\frac{b_{2k}}{2}+
c_{2n}c_{2t}3\Big(\frac{b_{2k}}{2}\Big)^2
$$
$$
\sim(b_{1k}+b_{2k})^2,\quad (f_k,f_r)=
$$
$$
\Big( \Big( x_{1n}x_{1t}\Big(D_{1k}^2+\frac{b_{1k}}{2}\Big)+
(x_{1n}x_{2t}+x_{2n}x_{1t})D_{1k}D_{2k}+x_{2n}x_{2t}D_{2k}^2
\Big){\bf 1},
$$
$$
 \Big(x_{1n}x_{1t}\Big(D_{1r}^2+\frac{b_{1r}}{2}\Big)+
(x_{1n}x_{2t}+x_{2n}x_{1t})D_{1r}D_{2r}+x_{2n}x_{2t}D_{2r}^2
\Big){\bf 1}\Big)=
$$
$$
c_{2n}c_{2t}\frac{b_{2k}}{2}\frac{b_{2r}}{2}\sim b_{2k}b_{2r}.
$$
Finally, we have
\begin{equation}
\label{(A_{r,s}),xx1,m=2}
(f_k,f_k) \sim (b_{1k}+b_{2k})^2, \quad(f_k,f_r)\sim b_{2k}b_{2r},\quad k\not=r.
\end{equation}
For $A_{(m)}=((f_k,f_r))_{k,r=1}^m,$ and
$b=(a_{11},a_{12},...,a_{1m})\in {\mathbb R}^m$ we have
\begin{equation}
\label{Gram-matr-x}
A_{(m)}= \gamma(f_1,f_2,...,f_m)=\left(
\begin{array}{cccc}
(f_1,f_1)&(f_1,f_2)&...&(f_1,f_m)\\
(f_2,f_1)&(f_2,f_2)&...&(f_2,f_m)\\
             &             &...&             \\
(f_m,f_1)&(f_m,f_2)&...&(f_m,f_m)
\end{array}
\right)=
\end{equation}
$$
\left(
\begin{array}{cccc}
(b_{11}+b_{21})^2& b_{21}b_{22}&...&b_{21}b_{2m}\\
 b_{22}b_{21}&(b_{12}+b_{22})^2&...&b_{22}b_{2m}\\
   &             &...&             \\
 b_{2m}b_{21} & b_{2m}b_{22}&...&(b_{1m}+b_{2m})^2
\end{array}
\right)=
$$
$$
\left(\!\!
\begin{array}{cccc}
b_{21}&0&...&0\\
0&b_{22}&...&0\\
 &      &...& \\
0&0&...&b_{2m}
\end{array}
\!\!\right)\!\! \left(\!\!
\begin{array}{cccc}
1+\lambda_1&1          &...&1\\
1          &1+\lambda_2&...&1\\
           &           &...& \\
1          &1&...&1+\lambda_m
\end{array}
\!\!\right)\!\! \left(\!\!
\begin{array}{cccc}
b_{21}&0&...&0\\
0&b_{22}&...&0\\
 &      &...& \\
0&0&...&b_{2m}
\end{array}
\!\!\right).
$$
At last, we have for $\mu={\rm diag}(b_{21},b_{22},...,b_{2m})^{-1}b=
(\frac{b_{11}}{b_{21}},\frac{b_{12}}{b_{22}},..,\frac{b_{1m}}{b_{2m}})$
$$
(A_{(m)}^{-1}b,b)=(D^{-1}(\lambda_1,\lambda_2,...,\lambda_m)\mu,\mu),\quad \lambda_k=\left(1+\frac{b_{1k}}{b_{2k}}\right)^2-1,\quad\mu_k=
\frac{b_{1k}}{b_{2k}}.
$$
 Using Lemma \ref{l7.7} for the operator $A_{2m+1}$, and the vector $b\in {\mathbb
R}^{2m+1}$ we obtain
$$
(A_{2m+1}^{-1}b,b)=\frac{\sum_{k=-m}^m
\frac{\left(\frac{b_{1k}}{b_{2k}}\right)^2}{\left(\frac{b_{1k}}{b_{2k}}+1\right)^2-1}
+ \sum_{-m\leq k<n\leq m}
 \frac {
\left(\frac{b_{1k}}{b_{2k}}-\frac{b_{1n}}{b_{2n}}\right)^2 } {
\left[\left(\frac{b_{1k}}{b_{2k}}+1\right)^2-1\right]
\left[\left(\frac{b_{1n}}{b_{2n}}+1\right)^2-1\right] }}
{\sum_{k=-m}^m\frac{1}{\left(\frac{b_{1k}}{b_{2k}}+1\right)^2-1}+1}=
$$
$\Delta(f^1_m,g^1_m)$ where

\begin{equation}
\label{f^1,g^1.1}
f_m^1=\Big(\frac{b_{1k}}{\sqrt{b_{1k}^2+2b_{1k}b_{2k}}}\Big)_{k=-m}^m,\quad g_m^1=\Big(\frac{b_{2k}}
{\sqrt{b_{1k}^2+2b_{1k}b_{2k}}}\Big)_{k=-m}^m,
\end{equation}
\qed\end{pf}
 The proof of  Lemma \ref{x2x2} is similar. We get $(A^{-1}b,b)=\Delta(f^2_m,g^2_m)$ where
\begin{equation}
\label{f^2,g^2.1}
f_m^2=\Big(\frac{b_{2k}}{\sqrt{b_{2k}^2+2b_{1k}b_{2k}}}\Big)_{k=-m}^m,\quad g_m^2=\Big(\frac{b_{1k}}{\sqrt{b_{2k}^2+2b_{1k}b_{2k}}}\Big)_
{k=-m}^m.
\end{equation}
%
%
{\it Acknowledgements.} {The author expresses his deep gratitude to
the Max Planck Institute for Mathematics for the financial grant and the hospitality in 2016-2017}.

\end{document}